\documentclass[12pt,reqno]{amsart}
\usepackage{amsthm,amsfonts,amssymb,euscript}

\newcommand{\bea}{\begin{eqnarray}}
\newcommand{\eea}{\end{eqnarray}}
\def\beaa{\begin{eqnarray*}}
\def\eeaa{\end{eqnarray*}}
\def\ba{\begin{array}}
\def\ea{\end{array}}
\def\be#1{\begin{equation} \label{#1}}
\def \eeq{\end{equation}}

\def\Lt{\widetilde{L}}

\def\a{{\alpha}}

\def\b{{\beta}}
\def\be{{\beta}}
\def\ga{\gamma}
\def\Ga{\Gamma}
\def\de{\delta}

\def\ep{\epsilon}
\def\eps{\epsilon}

\def\la{\lambda}

\def\si{\sigma}

\def\om{\omega}

\def\th{\theta}

\def\pr{{\partial}}
\def\al{\alpha}

\def\c{\cdot}

\def\MM{{\mathcal M}}
\def\NN{{\mathcal N}}

\def\HH{{\mathcal H}}

\def\OO{{\mathcal O}}

\def\NN{{\mathcal N}}

\def\KK{{\mathcal K}}
\def\Lie{{\mathcal L}}
\def\Lieh{\widehat{\Lie}}

\def\Lie{{\mathcal L}}

\def\D{{\bf D}}

\def\M{{\bf M}}

\def\R{{\bf R}}

\def\T{{\bf T}}

\def\g{{\bf g}}

\def\h{{\bf h}}
\def\piX{\,^{(X)}\pi}
\def\GaX{\,^{(X)}\Ga}
\def\piL{\,^{(L)}\pi}

\def\RRR{{\mathbb R}}

\def\f12{{\frac 1 2}}

\def\dual{{\,^*\,}}
\def\div{{\mbox div\,}}

\def\Lb{{\,\underline{L}}}

\def\trch{{\mbox tr}\, \chi}

\def\chih{{\hat \chi}}

\def\th{\theta}

\def\f{\widetilde{f}}

\def\um{\underline{u}}

\def\dual{\,^*}
\def\NNb{\underline{\NN}}
\def\ub{\underline{u}}

\newtheorem{theorem}{Theorem}[section]
\newtheorem{lemma}[theorem]{Lemma}
\newtheorem{proposition}[theorem]{Proposition}

\newtheorem{definition}[theorem]{Definition}
\newtheorem{remark}[theorem]{Remark}

\numberwithin{equation}{section}

\begin{document}
\title{On the local extension of Killing vector-fields in Ricci flat manifolds}
\author{Alexandru D. Ionescu}
\address{Princeton University}
\email{aionescu@math.princeton.edu}
\author{Sergiu Klainerman}
\address{Princeton University}
\email{seri@math.princeton.edu}

\thanks{The first author was supported in part by a Packard fellowship. The second author  was supported  in part by  NSF grant  0601186 as well as by  the 
 Fondation des Sciences Math\'ematiques  de Paris.}

\begin{abstract} We revisit the  extension problem for Killing vector-fields
in  smooth Ricci flat manifolds, and its  relevance to the  black hole rigidity problem. 
We  prove   both a   stronger version of the main local   extension  result  established  in \cite{AlIoKl},
  as well as two types of results  concerning non-extendibility.   In particular we      show 
     that one can find   local, stationary,   vacuum extensions of a Kerr solution $\KK(m,a)$, $0<a<m$, 
      in a  future neighborhood  of a point $p$ of the past  horizon, ($p$ not on  the bifurcation  sphere),   
        which admits no   extension of   the  Hawking vector-field   of $\KK(m,a)$.   This    result illustrates one  of the major   difficulties  one  faces in trying to extend Hawking's rigidity result to the more realistic setting of  smooth
     stationary  solutions of the Einstein vacuum equations; unlike in the analytic situation, one cannot hope to construct an additional symmetry of stationary solutions (as in Hawking's Rigidity Theorem) by relying only on  local information.
  
\end{abstract}
\maketitle
\tableofcontents

\section{Introduction}\label{intro}
In this paper we revisit the  extension problem for Killing vector-fields
in  smooth Ricci flat Lorentzian manifolds and its  relevance to the  black hole rigidity
 problem. In the most general
situation the problem can be stated as follows:

\textit{Assume $(\M,\g)$ is a given smooth pseudo-riemannian manifold, $O\subseteq\M$  is an open subset, and $Z$ is a  smooth Killing  vector-field in $O$.  Under what assumptions does $Z$ extend (uniquely) as a Killing vector-field in $\M$?}

A classical result\footnote{See \cite{No}. We rely here on the version of  the theorem given in \cite{Chrusc-1}.} of Nomizu establishes such a unique  extension   provided that the metric is real analytic, $\M$ and $O$ are connected and $\M$ is simply connected.  The result has been used, see \cite{HE} and \cite{Chrusc-1},    to reduce  the 
black hole rigidity problem, for  real analytic stationary solutions of the Einstein field equations, to the 
 simpler case  of axial symmetry  treated  by the Carter-Robinson theorem.  This reduction  
 has been    often   regarded as  decisive,  especially in the physics literature,    without  a clear understanding of the sweeping  simplification   power  of the analyticity assumption. 
  Indeed the remarkable thing about Nomizu's
 theorem,  to start with,  is the fact  the metric is not assumed to satisfy any specific equation. Moreover no assumptions are needed  about  the boundary of $O$
 in $\M$ and the result is  global with only minimal assumptions on the topology
 of $\M$ and $O$.
All these   are clearly wrong in the case of smooth manifolds $(\M, \g)$ which are not real analytic. To be able to say anything meaningful we need to  both  restrict the metric $\g$   by realistic  equations and   make specific assumptions about the boundary of $O$.   Local and global assumptions also need to be carefully separated.

\medskip
 In this paper we limit our  attention  to a purely local  description of the extension problem in the 
 smooth case. 
Throughout  the paper we  assume that $(\M, \g)$ is a non-degenerate Ricci flat, pseudo-riemannian metric  i.e.
\begin{equation}\label{Ricci}
\mathrm{\bf Ric}(\g)=0.
\end{equation}

We recall the following crucial concept.

 \begin{definition}\label{psconvexqual}
 A domain $O\subset \M$ is said to be strongly pseudo-convex at a boundary point  $p\in \pr O$
if it admits a strongly pseudo-convex defining function $f$  at $p$,  in the sense that there is an open neighborhood $U$ of $p$ in $\M$ and a smooth function $f:U\to\mathbb{R}$, $\nabla f(p)\neq 0$,  such that $O\cap U=\{x\in U: f(x)<0\}$ and
\begin{equation}\label{qual1.1}
\D^2f(X,X)(p)<0
\end{equation}
for any $X\neq 0\in T_p( \M)$ for which   $X(f)(p)=0$ and $\g_p(X,X)=0$.
\end{definition}

It is easy to see that this definition, in particular \eqref{qual1.1}, does not depend on the choice of the defining function $f$. The strong pseudo-convexity  condition is automatically satisfied  if the metric $\g$
  is Riemannian. It is also satisfied   for  Lorentzian  metrics $\g$
    if  $\pr O$ is space-like at $p$, but it imposes serious restrictions for time-like  hypersurfaces.    It  clearly  fails
 if $\pr O$ is null in a neighborhood of $p$.
  Indeed in that case  we can choose  the defining function
 $f$ to be optical,  i.e.,
 \begin{equation}\label{fail}
\D^\al f\D_\al f=0\quad
\end{equation}
 at all  points  of $\pr O$  in a neighborhood of  $p$, and thus, choosing $X^\al=\D^\al  f$, we  have,  $$X^\al X^\be\D_\al\D_\be f=\frac 1 2 X (\D^\al f\D_\al f)=0.$$
 
 Besides  a new extension result, see Theorem   \ref{extthm0} below,  the paper contains  two  local  counterexamples. In our main  such result, see  Theorem  \ref{mkathm},
  we show that   at any point  $p$  in the complement of the bifurcation sphere  of the  horizon of a Kerr spacetime $\KK(m,a), 0<a<m$, with  $\T, Z$ denoting  the usual  stationary and  axially symmetric Killing vector-fields of $\KK(m,a)$,  one can  find local extensions of  the Kerr metric, which coincide with   $\KK(m,a)$ inside the black hole, and such that only $\T$ extends as a  Killing  vector-field to  a full neighborhood of $p$. The condition $a>0$ is important in our proof, since our construction only works in the region where $\T$ is timelike, i.e. the ergo-region. It remains open whether a similar counterexample can be constructed  for the  Schwarzschild  spacetimes $\KK(m,0)$.
  
\bigskip

 We  first state   the following  extension theorem:

 \begin{theorem}\label{extthm0}
Assume that $(\M, \g) $ is a smooth $d$-dimensional  Ricci flat,  pseudo-rie\-mannian manifold  and  $O\subseteq \M$ is a strongly pseudo-convex domain at a point $p\in\pr O$. We assume that the metric $\g$ admits a smooth Killing vector-field $Z$ in $O$. Then $Z$ extends as a Killing vector-field for $\g$ to a neighborhood of the point $p$ in $\M$.
 \end{theorem}

Under more restrictive assumptions, a similar result was proved in \cite{AlIoKl}  as a key component of a theorem on the uniqueness of the Kerr solution in  \cite{AlIoKl2}.  In this  paper we  present  a different, more geometric proof,  which is valid in all dimensions and for all pseudo-riemannian metrics. More importantly, the proof we present here does not require that the vector-field $Z$ be tangent to the boundary $\partial O$ in a neighborhood of $p$, or the existence of a geodesic vector-field $L$, defined in a neighborhood of $p$, and commuting with $Z$ in $O$.

In applications, one would like to use Theorem \ref{extthm0} repeatedly and extend the Killing vector-field $Z$ to larger and larger open sets. For this it is important to understand the "size" of the implied neighborhood in the conclusion of the theorem, where the vector-field $Z$ extends. The proof shows that this neighborhood depends only on smoothness parameters of $\g$ and $f$ in a neighborhood of $p$ (see \eqref{quant2}), and a quantitative form of strong pseudo-convexity described in Lemma \ref{psconvexquan}. The neighborhood does not depend in any way on the vector-field $Z$ itself.

In view of Theorem \ref{extthm0}, Killing vector-fields extend locally across strongly pseudo-convex hypersurfaces in Ricci flat manifolds. A natural question is whether the strong pseudo-convexity condition is needed. We give a partial answer in Theorem \ref{nj2}: in general one cannot expect to extend a Killing vector-field across a null hypersurface in a $4$-dimensional Lorentz manifold.\footnote{Such a hypersurface is not strongly pseudo-convex, see the discussion before Theorem \ref{extthm0}.}

Our second main theorem provides a counterexample to extendibility, in  the setting of the black hole rigidity
problem.
 Let $(\mathcal{K}(m,a),\g)$ denote the (maximally extended) Kerr space-time of mass $m$ and angular momentum $ma$, $0\leq a<m$ (see \cite{HE} for definitions). Let $\M^{(end)}$ denote an asymptotic region, $\mathbf{E}=\mathcal{I}^-(\M^{(end)})\cap \mathcal{I}^+(\M^{(end)})$ the corresponding domain of outer communications, and $\mathcal{H}^-=\delta(\mathcal{I}^+(\M^{(end)})$ the boundary (event horizon) of the corresponding white hole\footnote{A similar statement can be made on  the future event horizon $\HH^+$.}. Let $\T=d/dt$ denote the stationary (timelike in $\M^{(end)}$) Killing vector-field  of $(\mathcal{K}(m,a),\g)$, and let $Z=d/d\phi$ denote its rotational (with closed orbits) Killing vector-field.

\begin{theorem}\label{mkathm}
Assume that $0<a<m$ and $U_0\subseteq \mathcal{K}(m,a)$ is an open set such that $$U_0\cap\mathcal{H}^-\cap\overline{\mathbf{E}}\neq\emptyset.$$ Then there is an open set $U\subseteq U_0$ diffeomorphic to the open unit ball $B_1\subseteq\mathbb{R}^4$, $U\cap \mathcal{H}^-\neq\emptyset$, and a smooth Lorentz metric $\widetilde{\g}$ in $U$ with the following properties:

(i)
\begin{equation}\label{mka1}
{}^{\widetilde{\g}}\mathbf{Ric}=0\,\,\text{ in }U,\qquad \Lie_\T{\widetilde{\g}}=0\,\,\text{ in }U,\qquad {\widetilde{\g}}=\g\,\,\text{ in }U\setminus\mathbf{E};
\end{equation}

(ii) the vector-field $Z=d/d\phi$ does not extend to a Killing vector-field for $\widetilde{\g}$, commuting with $\T$, in $U$.
\end{theorem}

In other words, one can modify the Kerr space-time smoothly, on one side of the horizon $\mathcal{H}^-$, in such a way that the resulting metric still satisfies the Einstein vacuum equations, has $\T=d/dt$ as a Killing vector-field, but does not admit an extension of the Killing vector-field $Z$. This    result illustrates one  of the major   difficulties  one  faces in trying to extend Hawking's rigidity result to the more realistic setting of  smooth
     stationary  solutions of the Einstein vacuum equations: unlike in the analytic situation, one cannot hope to construct an additional symmetry of stationary solutions of the Einstein-vacuum equations (as in Hawking's Rigidity Theorem) by relying only on the local information provided by the equations.\footnote{As mentioned earlier  a local version of Hawking's Rigidity Theorem was proved in \cite{AlIoKl}. The key additional information used in that paper is the existence of a regular bifurcation sphere, which is the smooth transversal intersection of two non-expanding horizons.}

     The rest of the paper is organized as follows: in section \ref{keyloc} we prove Theorem \ref{extthm0} and in section \ref{mka} we prove Theorem \ref{mkathm}. In section \ref{further} we consider extensions across null hypersurfaces in $4$-dimensional Lorentz manifolds and prove two more theorems: Theorem \ref{suff3}, which provides a criterion for extension of Killing vector-fields, and Theorem \ref{nj2}, which provides a general framework when extension is not possible.

\section{Proof of Theorem \ref{extthm0}}\label{keyloc} 
In   \cite{AlIoKl} and \cite{AlIoKl2} the extension of the Killing vector-field $Z$  was done according to the transport equation,
\bea
[L, Z]=c_0 L,\label{eq:commL}
\eea
where $\D_L L=0$ and $c_0$ constant. Consequently we  had to assume, in $O$,  that
$Z$ is not only Killing but that it also satisfies  the additional assumption 
\eqref{eq:commL} with respect to a geodesic non-vanishing vector-field  $L$. This could be arranged in the particular cases studied in \cite{AlIoKl} and \cite{AlIoKl2}, but imposes serious restrictions on $Z$ in the general case, particularly if $Z$ vanishes in a neighborhood of the point $p$. To avoid this restriction, in this paper we extend $Z$ according to the weaker condition 
\begin{equation}
\D_L \D_L Z=\R(L, Z) L\label{eq:extZ},
\end{equation}
which would follow easily from \eqref{eq:commL}, and is automatically satisfied  if $Z$ is Killing.
 \bigskip
 
More precisely, we construct first a smooth vector-field $L$ in a neighborhood of $p$ such that
\begin{equation*}
\D_LL=0,\qquad L(f)(p)=1,
\end{equation*}
and  extend  $Z$ to a neighborhood of $p$ by solving the second order  differential   system
 \eqref{eq:extZ}.
Therefore, after restricting to a small neighborhood of $p$, we may assume that $Z,L$ are smooth vector-fields in $\M$ with the properties
\begin{equation}\label{pr1.1}
\D_LL=0\,\,\text{ in }\M,\qquad L^\al L^\be(\D_\al\D_\be Z_\mu-Z^\rho\R_{\rho\al\be\mu})=0\,\,\text{ in }\M,\qquad \mathcal{L}_{Z}\g=0\,\,\text{ in } O.
\end{equation}
It remains to prove that the deformation tensor $\pi=\mathcal{L}_{Z}\g$ vanishes in a neighborhood of $p$. We   cannot do this however without establishing at the same time that the  tensor $\Lie_Z\R$ also vanishes identically in $\M$. Our strategy is to derive  a wave  equation for $\Lie_Z\R$, or rather a suitable modification of it, coupled with a number of transport equations along the integral curves of $L$ for various tensorial quantities including $\pi$ itself.  These equations will be used to prove that $\pi$ and $\Lie_Z\R$ have to vanish in a full neighborhood of $p$, provided that the strong pseudo-convexity assumption, which guarantees the unique continuation property, is satisfied.

\subsection{Tensorial equations} \label{sec:tens-eq}  We first consider the properties of $\Lie_Z \R$.  Observe that $\Lie_Z\R$
 verifies all  the algebraic  symmetries of $\R$ except the fact that,
  for an Einstein vacuum metric $\g$, $\R$ is traceless. We have instead,
  \beaa
  \g^{\a\ga} \Lie_Z\R_{\a\b\ga\de}&=&\pi^{\a\ga} \R_{\a\b\ga\de}.
  \eeaa
  To re-establish this property we can introduce (see also Chapter 7 in \cite{CKl}) modifications\footnote{Note however that, unlike   \cite{CKl},  our $B$ here is not symmetric.}  of $\Lie_Z \R$ of the form
  \beaa
  \Lieh_Z \R:=\Lie_Z\R-B\odot\R,
  \eeaa
  where, for any give $2$-tensor $B$,  we write,
  \beaa
 ( B\odot \R)_{\a\b\ga\de}:=B_\a\,^\la \R_{\la\b\ga\de}+B_\b\,^\la \R_{\a\la\ga\de}+
 B_\ga\,^\la \R_{\a\b\la\de}+B_\de\,^\la \R_{\a\b\ga\la}.
  \eeaa
  It is easy to check that, for any $2$-tensor $B$,   $B\odot \R$ verifies  all the algebraic
  symmetries of the general Riemann curvature tensor, i.e.
  \beaa
  &&(B\odot  \R)_{\a\b\ga\de}=-(B\odot  \R)_{\b\a\ga\de}=- (B\odot  \R)_{\a\b\de\ga}=(B\odot \R)_{\ga\de\a\b},\\
  && (B\odot \R)_{\a\b\ga\de}+  (B\odot \R)_{\a\ga\de\b}+  (B\odot  \R)_{\a\de\b\ga}=0.
  \eeaa
  Moreover, using the Einstein vacuum equations,
  \beaa
   \g^{\a\ga}(B\odot  \R)_{\a\b\ga\de}&=& B^{\mu\la}\big(\R_{\la\b\mu\de}+\R_{\mu\b\la\de}\big).
  \eeaa
  In particular  for any antisymmetric $B$,   $B\odot  \R$ is traceless, i.e. a Weyl field. We have proved the following:

\begin{proposition}\label{pro1}
Assume $\omega$ is an antisymmetric $2$-form in $\M$ and let
\begin{equation}\label{W:definition}
W:=\Lie_Z \R-\frac 1 2 (\pi+\omega)\odot \R.
\end{equation}
Then $W$ is a Weyl field in $\M$, i.e.
\begin{equation*}
\begin{split}
&W_{\a\b\ga\de}=-W_{\b\a\ga\de}=-W_{\a\b\de\ga}=W_{\ga\de\a\b},\\
&W_{\a\b\ga\de}+W_{\a\ga\de\b}+W_{\a\de\b\ga}=0,\\
&\g^{\al\ga}W_{\al\be\ga\de}=0.
\end{split}
\end{equation*}
\end{proposition}

We shall next establish a divergence equation for $W$. We do this by
commuting the  divergence equation for $\R$ with $\Lie_Z$. We rely on the following, see Lemma 7.1.3 in \cite{CKl}:

\begin{lemma}\label{le:commute}
For  arbitrary $k$-covariant  tensor-field  $V$ and vector-field  $X$ we have,
\begin{equation}\label{commute}
\D_\b(\mathcal{L}_X V_{\al_1\ldots\al_k})-\mathcal{L}_X(\D_\b V_{\al_1\ldots\al_k})=\sum_{j=1}^k\GaX_{\al_j\b \rho}V_{\al_1\ldots\,\,\,\,\ldots\al_k}^{\,\,\,\,\,\,\,\,\,\,\,\rho},
\end{equation}
where $\piX=\Lie_X\g$  is the deformation tensor of $X$ and,
\begin{equation*}
\GaX_{\al\be\mu}=\frac 1 2 (\D_\al\piX_{\be\mu}+\D_\be\piX_{\al\mu}-\D_\mu\piX_{\al\be}).
\end{equation*}
\end{lemma}

\begin{definition}\label{pro0}
We denote   $\pi=\,^{(Z)}\pi$ and $\Ga=\, ^{(Z)}\Ga$ the corresponding
tensors associated to the vector-field $Z$. We  also
 denote  $\piL=H$. We also introduce the tensors,
\begin{equation*}
\begin{split}
&P_{\a\b\mu}=\D_\a\pi_{\b\mu}-\D_\b\pi_{\a\mu}-\D_\mu \om_{\a\b},\\
&B_{\al\be}=\frac 1 2 (\pi_{\al\be}+\omega_{\al\be}),\\
&\dot{B}_{\al\be}=L^\rho\D_\rho B_{\al\be},\\
&W_{\al\be\ga\de}=(\Lie_Z\R)_{\al\be\ga\de}-(B\odot\R)_{\al\be\ga\de}.
\end{split}
 \end{equation*}
All these tensors depend on the $2$-form $\omega$, which will be defined later (see \eqref{maincor}) to achieve a key cancellation in the proof of the transport equation \eqref{yi3}.
 \end{definition}

Using Lemma \ref{le:commute} we can now prove the following:

\begin{lemma}\label{pro2}
The Weyl field $W$ verifies the divergence equation
\begin{equation}\label{eq:J}
\begin{split}
\D^\a W_{\a\b\ga\de}=\frac{1}{2}\big(&B^{\mu\nu}\D_\nu \R_{\mu\b\ga\de}+\g^{\mu\nu}P_{\mu\rho\nu}\R^\rho\,_{\b\ga\de}\\
  &+P_{\b\nu\mu}\R^{\mu\nu}\,_{\ga\de}+P_{\ga\nu\mu}   \R^{\mu}\,_\b\,^\nu\,_{\de}
  +P_{\de\nu\mu}\R^{\mu}\,_{\b\ga}\,^\nu\big).
\end{split}
\end{equation}
\end{lemma}

\begin{proof}[Proof of Lemma \ref{pro2}] Using Lemma \ref{le:commute} and the identity $\D^\al\R_{\al\be\ga\de}=0$ (which is a consequence of the Einstein vacuum equations), we easily deduce
\begin{equation*}
\begin{split}
\D^\al\mathcal{L}_Z\R_{\al\be\ga\de}&=\g^{\al\mu}\D_\mu\mathcal{L}_Z\R_{\al\be\ga\de}\\
&=\g^{\al\mu}\big(\mathcal{L}_Z\D_\mu\R_{\al\be\ga\de}+\Gamma_{\al\mu\rho}{\R^\rho}_{\be\ga\de}+\Gamma_{\be\mu\rho}{{\R_\al}^\rho}_{\ga\de}+\Gamma_{\ga\mu\rho}{{\R_{\al\be}}^\rho}_\de+\Gamma_{\de\mu\rho}{\R_{\al\be\ga}}^\rho\big)\\
&=\pi^{\al\mu}\D_\mu\R_{\al\be\ga\de}+{\Gamma^\mu}_{\mu\rho}{\R^\rho}_{\be\ga\de}+\Gamma_{\be\mu\rho}{\R^{\mu\rho}}_{\ga\de}+\Gamma_{\ga\mu\rho}{{{\R^\mu}_{\be}}^\rho}_\de+\Gamma_{\de\mu\rho}{{\R^\mu}_{\be\ga}}^\rho.
\end{split}
\end{equation*}
Using the definition and the Einstein vacuum equations, we derive
  \begin{equation*}
  \begin{split}
  &\D^\al(B\odot \R)_{\a\b\ga\de}\\
  &=B^{\al\la}\D_\al\R_{\la\be\ga\de}+\D^\al B_{\a\la} {\R^\la}_{\b\ga\de}+\D_\al B_{\b\la} {\R^{\al\la}}_{\ga\de}+
 \D_\al B_{\ga\la} {{{\R^\al}_{\b}}^{\la}}_\de+\D_\al B_{\de\la} {{\R^\al}_{\b\ga}}^\la,
 \end{split}
  \end{equation*}
  for any $2$-tensor $B$. Thus, if $B=(1/2)(\pi+\omega)$,
  \begin{equation*}
  \begin{split}
  \D^\al W_{\al\be\ga\de}&=(\pi^{\mu\nu}-B^{\mu\nu})\D_\mu\R_{\nu\be\ga\de}+\g^{\mu\nu}(\Gamma_{\mu\nu\rho}-\D_\nu B_{\mu\rho}){\R^\rho}_{\b\ga\de}\\
  &+(\Gamma_{\be\mu\nu}-\D_\mu B_{\be\nu}){\R^{\mu\nu}}_{\ga\delta}+(\Gamma_{\ga\mu\nu}-\D_\mu B_{\ga\nu}){{{\R^\mu}_{\be}}^\nu}_\de+(\Gamma_{\de\mu\nu}-\D_\mu B_{\de\nu}){{\R^\mu}_{\be\ga}}^\nu.
  \end{split}
  \end{equation*}
  We observe now that
  \begin{equation*}
  \Gamma_{bac}-\D_a B_{bc}=\frac 1 2 P_{bca},
  \end{equation*}
  which completes the proof of the lemma.
\end{proof}

We now look for transport equations for the tensor-fields $B,P$ appearing in \eqref{eq:J}, of the form,
\beaa
\D_L(B, P)=\MM(W, B, P),
\eeaa
with  the notation $\MM(W, B, P)$ explained below.

 \begin{definition}
 By convention, we let $\mathcal{M}({}^{(1)}B,\ldots,{}^{(k)}B)$ denote any smooth ``multiple'' of the tensors ${}^{(1)}B,\ldots,{}^{(k)}B$, i.e. any tensor of the form
\begin{equation}\label{mnotation}
\mathcal{M}({}^{(1)}B,\ldots,{}^{(k)}B)_{\al_1\ldots\al_r}={}^{(1)}B_{\be_1\ldots\be_{m_1}}{}^{(1)}{C_{\al_1\ldots\al_r}}^{\be_1\ldots\be_{m_1}}+\ldots +{}^{(k)}B_{\be_1\ldots\be_{m_k}}{}^{(k)}{C_{\al_1\ldots\al_r}}^{\be_1\ldots\be_{m_k}},
\end{equation}
for some smooth tensors ${}^{(1)}C,\ldots, {}^{(k)}C$ in $\M$.
 \end{definition}
 It turns out in fact that we need to include also a transport equation
 for $\dot{B}=\D_L B$. Thus we look for equations of the type,
 \beaa
\D_L(B,\dot{B}, P)=\MM(W, B, \dot{B}, P).
\eeaa
  We start with a lemma.

  \begin{lemma}\label{pro3}
  Given the vector-field $Z$, extended to $\M$ by \eqref{pr1.1},
  we have
  \begin{equation}\label{Lpi}
  L^\be\pi_{\al\be}=0\qquad\text{ in }\M.
  \end{equation}
  Moreover, if we define $\omega$ in $\M$ as the  solution of the transport equation
  \begin{equation}\label{maincor}
\D_L\omega_{\al\be}=\pi_{\al\rho}\D_\be L^\rho-\pi_{\be\rho}\D_\al L^\rho,
\end{equation}
with $\omega=0$ in $O$, then
\begin{equation}\label{LP}
L^\mu P_{\al\be\mu}=0,\qquad L^\be\om_{\al\be}=0\qquad\text{ in \,\, } \M.
\end{equation}
  \end{lemma}

  \begin{proof}[Proof of Lemma \ref{pro3}]
  We first remark that $L^\al L^\be\D_\al Z_\be=0$ in $\M$. Indeed, using \eqref{pr1.1},
  \begin{equation*}
 L^\rho\D_\rho(L^\a L^\b \D_\a Z_\b)=L^\rho L^\a L^\b \D_\rho\D_\a Z_\b=L^\be L^\al L^\rho Z^\nu\R_{\nu\al\rho\be}=0.
  \end{equation*}
 Since $L^\a L^\b \D_\a Z_\b=0$ in $O$ we deduce that
  \begin{equation}\label{step1}
  L^\al L^\be\D_\al Z_\be=0\qquad \mbox{in} \,\, \M.
  \end{equation}

We prove now \eqref{Lpi}. Using \eqref{pr1.1} and \eqref{step1} we compute
\begin{equation*}
\begin{split}
L^\rho\D_\rho(L^\be\pi_{\al\be})&=L^\rho L^\be(\D_\rho\D_\be Z_\al+\D_\rho\D_\al Z_\be)\\
&=L^\rho L^\be Z^\mu\R_{\mu\rho\be\al}+L^\rho L^\be\D_\al\D_\rho Z_\be+L^\rho L^\be Z^\mu\R_{\rho\al\be\mu}\\
&=\D_\al(L^\rho L^\be\D_\rho Z_\be)-L^\be\D_\rho Z_\be\D_\al L^\rho-L^\rho\D_\rho Z_\be\D_\al L^\be\\
&=-L^\be\pi_{\mu\be}\D_\al L^\mu.
\end{split}
\end{equation*}
Since $L^\be\pi_{\al\be}$ vanishes in $O$, it follows that $L^\be\pi_{\al\be}$ vanishes in $\M$, as desired.

The first identity in \eqref{LP} follows from the definitions of $\omega$ and $P$ and the identity \eqref{Lpi}:
\begin{equation*}
L^\mu P_{\al\be\mu}=L^\mu \D_\al\pi_{\be\mu}-L^\mu\D_\be\pi_{\al\mu}-L^\mu\D_\mu\om_{\al\be}\\
=-\pi_{\be\mu}\D_\al L^\mu+\pi_{\al\mu}\D_\be L^\mu-L^\mu\D_\mu\om_{\al\be}=0.
\end{equation*}
To prove the second identity, we compute, using the definition \eqref{maincor} and the identities $L^\be\pi_{\be\rho}=0$ and $\D_LL=0$,
\begin{equation*}
\D_L(L^\be\omega_{\al\be})=L^\rho L^\be\D_\rho\omega_{\al\be}=L^\be(\pi_{\al\rho}\D_\be L^\rho-\pi_{\be\rho}\D_\al L^\rho)=0.
\end{equation*}
Since $L^\be\omega_{\al\be}$ vanishes in $O$, it follows that $L^\be\omega_{\al\be}$ vanishes in $\M$, as desired.
\end{proof}

We derive now our main transport equations for the tensors $\dot{B}$ and $P$.

\begin{lemma}\label{pro5}
In $\M$ we have
\begin{equation}\label{yi2}
   \D_L\dot{B}_{\a\b}=L^\mu L^\nu   (\Lie_Z \R)_{\mu\al\be\nu}-2\dot{B}_{\nu\be}\D_\al L^\nu-{\pi_\be}^\rho L^\mu L^\nu\R_{\mu\al\rho\nu}
   \end{equation}
and
\begin{equation}\label{yi3}
\D_L P_{\al\be\mu}=2L^\nu W_{\al\be\mu\nu}+2L^\nu {B_{\mu}}^\rho\R_{\al\be\rho\nu}-\D_\mu L^\rho P_{\al\be\rho}.
\end{equation}
  \end{lemma}

  \begin{proof}[Proof of Lemma \ref{pro5}] We have
  \begin{equation}\label{do1}
  \D_L\dot{B}_{\al\be}=L^\mu L^\nu\D_\mu\D_\nu B_{\al\be}=(1/2)[L^\mu L^\nu\D_\mu\D_\nu \pi_{\al\be}+L^\mu \D_\mu(L^\nu\D_\nu \om_{\al\be})].
  \end{equation}
  We calculate
  \begin{equation*}
  \begin{split}
  L^\mu L^\nu\D_\mu\D_\nu \D_\al Z_\be&=L^\mu L^\nu\D_\mu(\D_\al \D_\nu Z_\be+\R_{\nu\al\be\rho}Z^\rho)\\
  &=L^\mu L^\nu(\D_\al\D_\mu \D_\nu Z_\be+{\R_{\mu\al\nu}}^\rho\D_\rho Z_\be+{\R_{\mu\al\be}}^\rho\D_\nu Z_\rho)\\
  &+L^\mu L^\nu Z^\rho\D_\mu\R_{\nu\al\be\rho}+L^\mu L^\nu \D_\mu Z_\rho{\R_{\nu\al\be}}^\rho.
  \end{split}
  \end{equation*}
  Using \eqref{pr1.1} and the general identity
  \begin{equation}\label{do6}
  \D_a\D_b Z_c=\R_{cbad}Z^d+\Gamma_{abc},
  \end{equation}
  we calculate
  \begin{equation*}
  \begin{split}
  L^\mu L^\nu\D_\al\D_\mu \D_\nu Z_\be&=\D_\al(L^\mu L^\nu Z^\rho \R_{\rho\mu\nu\be})-\D_\al(L^\mu L^\nu)(\R_{\be\nu\mu\rho}Z^\rho +\Gamma_{\nu\mu\be})\\
  &=L^\mu L^\nu Z^\rho \D_\al\R_{\rho\mu\nu\be}+L^\mu L^\nu \D_\al Z_\rho {\R^\rho}_{\mu\nu\be}-\Gamma_{\nu\mu\be}\D_\al(L^\mu L^\nu).
  \end{split}
  \end{equation*}
  Thus
  \begin{equation*}
  \begin{split}
  L^\mu L^\nu\D_\mu\D_\nu \D_\al Z_\be&=L^\mu L^\nu Z^\rho(\D_\mu\R_{\nu\al\be\rho}+\D_\al\R_{\rho\mu\nu\be})-\Gamma_{\nu\mu\be}\D_\al(L^\mu L^\nu)\\
  &+L^\mu L^\nu(\D_\al Z^\rho \R_{\mu\rho\be\nu}+2\D_\nu Z^\rho \R_{\mu\al\be\rho}+\D_\be Z^\rho \R_{\mu\al\rho\nu})+L^\mu L^\nu{\R_{\mu\al\nu}}^\rho\pi_{\rho\be}.
  \end{split}
  \end{equation*}
  Since
  \begin{equation*}
  (\D_\mu\R_{\nu\al\be\rho}+\D_\al\R_{\rho\mu\nu\be})+(\D_\mu\R_{\nu\be\al\rho}+\D_\be\R_{\rho\mu\nu\al})=\D_\rho\R_{\mu\be\al\nu}+\D_\rho\R_{\mu\al\be\nu},
  \end{equation*}
  it follows that
  \begin{equation*}
  \begin{split}
  L^\mu L^\nu\D_\mu\D_\nu \pi_{\al\be}&=L^\mu L^\nu[(\Lie _Z\R)_{\mu\al\be\nu}+(\Lie _Z\R)_{\mu\be\al\nu}]\\
  &+L^\mu L^\nu({\R_{\mu\al\nu\rho}}{\pi_{\be}}^\rho+{\R_{\mu\be\nu\rho}}{\pi_{\al}}^\rho)-[\Gamma_{\nu\mu\be}\D_\al(L^\mu L^\nu)+\Gamma_{\nu\mu\al}\D_\be(L^\mu L^\nu)].
  \end{split}
  \end{equation*}
 Using the identity \eqref{Lpi} and the definitions we calculate
  \begin{equation*}
  \Gamma_{\nu\mu\be}\D_\al(L^\mu L^\nu)=L^\nu\D_\nu\pi_{\mu\be}\D_\al L^\mu-\D_\al L^\mu\D_\mu L^\nu\pi_{\nu\be}+\D_\al L^\mu\D_\be L^\nu\pi_{\mu\nu}=2\D_\al L^\mu\dot{B}_{\mu\be}.
  \end{equation*}
  Therefore
  \begin{equation}\label{do2}
  \begin{split}
   L^\mu L^\nu\D_\mu\D_\nu \pi_{\al\be}&=2L^\mu L^\nu(\Lie _Z\R)_{\mu\al\be\nu}-2\D_\al L^\mu\dot{B}_{\mu\be}-2\D_\be L^\mu\dot{B}_{\mu\al}\\
   &+L^\mu L^\nu({\R_{\mu\al\nu\rho}}{\pi_{\be}}^\rho+{\R_{\mu\be\nu\rho}}{\pi_{\al}}^\rho).
   \end{split}
  \end{equation}

  Using again \eqref{Lpi} and the definitions we calculate
  \begin{equation}\label{do3}
  \begin{split}
  L^\mu \D_\mu&(L^\nu\D_\nu \om_{\al\be})=L^\mu \D_\mu(\pi_{\al\nu}\D_\be L^\nu-\pi_{\be\nu}\D_\al L^\nu)\\
  &=(L^\mu \D_\mu\pi_{\nu\al}\D_\be L^\nu+{\pi_\al}^\rho L^\mu \D_\mu\D_\be L_\rho)- (L^\mu \D_\mu\pi_{\nu\be}\D_\al L^\nu+{\pi_\be}^\rho L^\mu \D_\mu\D_\al L_\rho)\\
  &=(2\D_\be L^\nu\dot{B}_{\nu\al}+{\pi_\al}^\rho L^\mu L^\nu\R_{\mu\be\rho\nu})-(2\D_\al L^\nu\dot{B}_{\nu\be}+{\pi_\be}^\rho L^\mu L^\nu\R_{\mu\al\rho\nu}).
  \end{split}
  \end{equation}
  The desired identity \eqref{yi2} follows from \eqref{do1}, \eqref{do2}, and \eqref{do3}.

  We prove now \eqref{yi3}. It follows from Definition \ref{pro0} that
  \begin{equation}\label{do5}
  \begin{split}
  L^\rho\D_\rho P_{\al\be\mu}&=L^\rho\D_\rho\D_\al(\D_\be Z_\mu+\D_\mu Z_\be)-L^\rho\D_\rho\D_\be(\D_\al Z_\mu+\D_\mu Z_\al)-L^\rho\D_\rho\D_\mu\omega_{\al\be}\\
  &=L^\rho\D_\rho(\R_{\al\be\mu\nu}Z^\nu)+(L^\rho\D_\rho\D_\al\D_\mu Z_\be-L^\rho\D_\rho\D_\be\D_\mu Z_\al)-L^\rho\D_\rho\D_\mu\omega_{\al\be}.
  \end{split}
  \end{equation}
  We calculate as before, using \eqref{do6}
  \begin{equation*}
  \begin{split}
  L^\rho&\D_\rho\D_\al\D_\mu Z_\be=L^\rho\D_\rho(\R_{\al\mu\be\nu}Z^\nu)+L^\rho\D_\rho\D_\mu\D_\al Z_\be\\
  &=L^\rho\D_\rho(\R_{\al\mu\be\nu}Z^\nu)+L^\rho{\R_{\rho\mu\al}}^\nu\D_\nu Z_\be+L^\rho{\R_{\rho\mu\be}}^\nu\D_\al Z_\nu+L^\rho\D_\mu\D_\rho\D_\al Z_\be\\
  &=L^\rho\D_\rho(\R_{\al\mu\be\nu}Z^\nu)+L^\rho{\R_{\rho\mu\al}}^\nu\D_\nu Z_\be+L^\rho{\R_{\rho\mu\be}}^\nu\D_\al Z_\nu+L^\rho\D_\mu(\R_{\be\al\rho\nu}Z^\nu+\Gamma_{\rho\al\be}).
  \end{split}
  \end{equation*}
  Using \eqref{Lpi} we calculate
  \begin{equation*}
  \begin{split}
  2L^\rho\D_\mu\Gamma_{\rho\al\be}&=L^\rho\D_\mu(\D_\rho\pi_{\al\be}+\D_\al\pi_{\rho\be}-\D_\be\pi_{\al\rho})\\
  &=L^\rho\D_\mu\D_\rho\pi_{\al\be}+\D_\mu L^\rho(\D_\be\pi_{\rho\al}-\D_\al\pi_{\rho\be})+\D_\mu(\D_\be L^\rho\pi_{\rho\al}-\D_\al L^\rho\pi_{\rho\be}).
  \end{split}
  \end{equation*}
  The last two identities and the definitions show that
  \begin{equation*}
  \begin{split}
  L^\rho&\D_\rho(\R_{\al\be\mu\nu}Z^\nu)+(L^\rho\D_\rho\D_\al\D_\mu Z_\be-L^\rho\D_\rho\D_\be\D_\mu Z_\al)\\
  &=2L^\rho\D_\rho(\R_{\al\be\mu\nu}Z^\nu)+2L^\rho\D_\mu(\R_{\be\al\rho\nu}Z^\nu)\\
  &+L^\rho{\R_{\rho\mu\al}}^\nu(\D_\nu Z_\be-\D_\be Z_\nu)-L^\rho{\R_{\rho\mu\be}}^\nu(\D_\nu Z_\al-\D_\al Z_\nu)\\
  &+\D_\mu L^\rho(\D_\be\pi_{\rho\al}-\D_\al\pi_{\rho\be})+\D_\mu(\D_\be L^\rho\pi_{\rho\al}-\D_\al L^\rho\pi_{\rho\be})\\
  &=2L^\rho (Z^\nu\D_\nu\R_{\al\be\mu\rho}+\D_\rho Z^\nu\R_{\al\be\mu\nu}+\D_\mu Z^\nu\R_{\al\be\nu\rho}+\D_\al Z^\nu\R_{\nu\be\mu\rho}+\D_\be Z^\nu\R_{\al\nu\mu\rho})\\
  &+L^\rho{\R_{\rho\mu\al}}^\nu\pi_{\be\nu}-L^\rho{\R_{\rho\mu\be}}^\nu\pi_{\al\nu}-\D_\mu L^\rho (P_{\al\be\rho}+\D_\rho\omega_{\al\be})+\D_\mu(L^\rho\D_\rho\omega_{\al\be})\\
  &=2L^\rho(\mathcal{L}_Z\R)_{\al\be\mu\rho}-L^\rho{\pi_\al}^\nu\R_{\nu\be\mu\rho}-L^\rho{\pi_{\be}}^\nu\R_{\al\nu\mu\rho}-\D_\mu L^\rho P_{\al\be\rho}+L^\rho\D_\mu\D_\rho\omega_{\al\be}.
  \end{split}
  \end{equation*}
  Using \eqref{do3} it follows that
  \begin{equation*}
  L^\rho\D_\rho P_{\al\be\mu}+\D_\mu L^\rho P_{\al\be\rho}=2L^\rho(\mathcal{L}_Z\R)_{\al\be\mu\rho}-2L^\rho{B_\al}^\nu\R_{\nu\be\mu\rho}-2L^\rho{B_{\be}}^\nu\R_{\al\nu\mu\rho},
  \end{equation*}
  which is equivalent to \eqref{yi3} (since $L^\rho B_{\rho\nu}=0$, see Lemma \ref{pro3}).
  \end{proof}

Finally, we derive a wave equation for the tensor $W$.

\begin{lemma}\label{pro7}
With the notation in \eqref{mnotation},
\begin{equation*}
\D^\rho\D_\rho W_{\al\be\mu\nu}=\mathcal{M}(B,\D B,P,\D P,W)_{\al\be\mu\nu}.
\end{equation*}
\end{lemma}

\begin{proof}[Proof of Lemma \ref{pro7}] We use the identity
\begin{equation}\label{yi30}
\begin{split}
\D^\si\D_\si\R_{\al_1\al_2\al_3\al_4}&=\R_{\si\rho\al_3\al_4}{{\R^\si}_{\al_1\al_2}}^\rho+\R_{\si\al_2\rho\al_4}{{\R^\si}_{\al_1\al_3}}^\rho+\R_{\si\al_2\al_3\rho}{{\R^\si}_{\al_1\al_4}}^\rho\\
&-\R_{\si\rho\al_3\al_4}{{\R^\si}_{\al_2\al_1}}^\rho-\R_{\si\al_1\rho\al_4}{{\R^\si}_{\al_2\al_3}}^\rho-\R_{\si\al_1\al_3\rho}{{\R^\si}_{\al_2\al_4}}^\rho,
\end{split}
\end{equation}
which is a well-known consequence of the Einstein vacuum equations. Using Lemma \ref{le:commute},
\begin{equation*}
\D_\si (\Lie_Z\R)_{\al_1\al_2\al_3\al_4}=\mathcal{L}_Z(\D_\si\R_{\al_1\al_2\al_3\al_4})+\sum_{j=1}^4\Gamma_{\al_j\si\rho}\R_{\al_1\ldots\,\,\,\,\ldots\al_4}^{\,\,\,\,\,\,\,\,\,\,\,\,\,\rho},
\end{equation*}
and then
\begin{equation*}
\square_\g (\Lie_Z\R)_{\al_1\al_2\al_3\al_4}=\D^\si\mathcal{L}_Z(\D_\si\R_{\al_1\al_2\al_3\al_4})+\sum_{j=1}^4[(\D^\si\Gamma_{\al_j\si\rho})\R_{\al_1\ldots\,\,\,\,\ldots\al_4}^{\,\,\,\,\,\,\,\,\,\,\,\,\,\rho}+\Gamma_{\al_j\si\rho}\D^\si\R_{\al_1\ldots\,\,\,\,\ldots\al_4}^{\,\,\,\,\,\,\,\,\,\,\,\,\,\rho}].
\end{equation*}
Therefore, after using Lemma \ref{le:commute} to commute derivatives again and \eqref{yi30}, the equation for $\square_\g (\Lie_Z\R)$ above can be written, in schematic notation,
\begin{equation*}
\square_\g (\Lie_Z\R)_{\al_1\ldots\al_4}=\sum_{j=1}^4[(\D^\si\Gamma_{\al_j\si\rho})\R_{\al_1\ldots\,\,\,\,\ldots\al_4}^{\,\,\,\,\,\,\,\,\,\,\,\,\,\rho}]+\mathcal{M}(\pi,\D\pi,(\Lie_Z\R))_{\al_1\ldots\al_4}.
\end{equation*}
In view of the definition,
\begin{equation*}
\square_\g(B\odot\R)_{\al_1,\ldots,\al_4}=\sum_{j=1}^4\D^\si\D_\si B_{\al_j\rho}\R_{\al_1\ldots\,\,\,\,\ldots\al_4}^{\,\,\,\,\,\,\,\,\,\,\,\,\,\rho}+\mathcal{M}(B,\D B)_{\al_1\ldots\al_4}.
\end{equation*}
Using also $\pi=\mathcal{M}(B)$, $\D\pi=\mathcal{M}(\D B)$, and $\Lie_Z \R=\mathcal{M}(B,W)$, it follows that
\begin{equation*}
\square_\g W_{\al_1\ldots\al_4}=\sum_{j=1}^4[\D^\si(\Gamma_{\al_j\si\rho}-\D_\si B_{\al_j\rho})\R_{\al_1\ldots\,\,\,\,\ldots\al_4}^{\,\,\,\,\,\,\,\,\,\,\,\,\,\rho}]+\mathcal{M}(B,\D B,W)_{\al_1\ldots\al_4}.
\end{equation*}
The lemma follows using the identity $\Gamma_{\al\be\mu}-\D_\be B_{\al\mu}=(1/2)P_{\al\mu\be}$.
\end{proof}

We summarize some of the main results in this subsection in the following proposition:

\begin{proposition}\label{equ1}
We assume that $O\subseteq \M,L,Z$ are as defined at the beginning of this section, and satisfy \eqref{pr1.1}. In $\M$ we define
\begin{equation*}
\pi_{\al\be}=\D_\al Z_\be+\D_\be Z_\al.
\end{equation*}
We define the smooth antisymmetric  tensor $\omega_{\al\be}$ in $\M$ as the solution of the equation
\begin{equation*}
\D_L\omega_{\al\be}=\pi_{\al\rho}\D_\be L^\rho-\pi_{\be\rho}\D_\al L^\rho,\qquad\omega=0\text{ in }O.
\end{equation*}
We also  define the smooth tensors
\begin{equation*}
\begin{split}
&P_{\a\b\mu}=\D_\a\pi_{\b\mu}-\D_\b\pi_{\a\mu}-\D_\mu \om_{\a\b},\\
&B_{\al\be}=\frac 1 2 (\pi_{\al\be}+\omega_{\al\be}),\\
&\dot{B}_{\al\be}=L^\rho\D_\rho B_{\al\be},\\
&W_{\al\be\ga\de}=(\Lie_Z\R)_{\al\be\ga\de}-(B\odot\R)_{\al\be\ga\de}.
\end{split}
\end{equation*}
Then the following equations hold in $\M$:
\begin{equation}\label{schembig}
\begin{split}
&\D^\al W_{\al\be\ga\de}=\mathcal{M}(B,\dot{B},P,W)_{\be\ga\de},\\
&\D_LB=\mathcal{M}(B,\dot{B},P,W),\quad\D_L\dot{B}=\mathcal{M}(B,\dot{B},P,W),\quad\D_LP=\mathcal{M}(B,\dot{B},P,W),\\
&\square W=\mathcal{M}(B,\D B,\dot{B},\D\dot{B},P,\D P,W,\D W).
\end{split}
\end{equation}
where $\mathcal{M}({}^{(1)}B,\ldots,{}^{(k)}B)$ is defined as in \eqref{mnotation}.
\end{proposition}

\subsection{Carleman inequalities and the local extension theorem}

Motivated by the identities summarized in Proposition \ref{equ1}, we consider solutions of systems of equations of the form
\begin{equation*}
\begin{cases}
&\square_\g S=\mathcal{M}({}^{(1)}B,\ldots,{}^{(k)}B,S,\D S)\\
&\D_L{}^{(i)}B=\mathcal{M}({}^{(1)}B,\ldots,{}^{(k)}B,S,\D S),\qquad i=1,\ldots,k.
\end{cases}
\end{equation*}
We would like to prove that a solution $S,{}^{(1)}B,\ldots,{}^{(k)}B$ of such a system which vanishes on one side of a suitable hypersurface has to vanish in a neighborhood of the hypersurface. Such a result depends, of course, on convexity and non-degeneracy properties of the hypersurface. We recall,
 see definition \ref{psconvexqual},  that a domain $O$  is strongly pseudo-convex at  a  boundary point  $p$ if  there exists  a defining function $f$ at $p$, $df(p)\neq 0$  which verifies,
\begin{equation}\label{qual1}
\D^2f(X,X)(p)<0\text{ if }X\neq 0\in T_p(\M)\,\, \text{ satisfies }\,\, \g_p(X,X)=X(f)(p)=0.
\end{equation}

We are now ready to prove Theorem \ref{extthm0}. We use the covariant equations derived in Proposition \ref{equ1} (see \eqref{schembig}) and Carleman inequalities. We introduce a smooth system of coordinates $\Phi^p=(x^1,\ldots,x^d):B_1\to B_1(p)$, $\Phi^p(0)=p$, where $B_r=\{x\in\mathbb{R}^d:|x|<r\}$, $r>0$, and $B_1(p)$ is an open neighborhood of $p$ in $\M$. Let $\partial_1,\ldots,\partial_d$ denote the induced coordinate vector-fields in $B_1(p)$ and let $B_r(p)=\Phi^p(B_r)$, $r\in(0,1]$. For any smooth function $\phi:B\to\mathbb{C}$, where $B\subseteq B_1(p)$ is an open set, and $j=0,1,\ldots$, we define
\begin{equation}\label{quant0}
|\partial^j\phi(x)|=\sum_{\al_1,\ldots,\al_j=1}^d|\partial_{\al_1}\ldots\partial_{\al_j}\phi(x)|,\qquad x\in B.
\end{equation}
We assume that
\begin{equation}\label{quant1}
\g_{\al\be}(p)=\mathrm{diag}(-1,\ldots,-1,1,\ldots,1).
\end{equation}
We assume also that, for some constant $A\geq 1$,
\begin{equation}\label{quant2}
\sup_{x\in B_1(p)}\sum_{j=1}^6\sum_{\al,\be=1}^d|\partial^j\g_{\al\be}(x)|+\sup_{x\in B_1(p)}\sum_{j=1}^4|\partial^j f(x)|\leq A.
\end{equation}

We use the system of coordinates $\Phi^p$ in the neighborhood of the point $p$, and evaluate all the tensor-fields in the frame of coordinate vector-fields $\partial_1,\ldots\partial_d$. In view of the equations \eqref{schembig}, for Theorem \ref{extthm0} it suffices to prove the following:

\begin{lemma}\label{extendedCarl2}
Assume that $\delta_0>0$ and $G_i,H_j:B_{\delta_0}(p)\to\mathbb{C}$ are smooth functions, $i=1,\ldots,I$, $j=1,\ldots,J$, that satisfies the differential inequalities
\begin{equation}\label{ext68}
\begin{cases}
&|\square_\g G_i|\leq M\sum_{l=1}^I(|G_l|+|\partial^1G_l|)+M\sum_{m=1}^J|H_m|;\\
&|L(H_j)|\leq M\sum_{l=1}^I(|G_l|+|\partial^1G_l|)+M\sum_{m=1}^J|H_m|,
\end{cases}
\end{equation}
for any $i=1,\ldots,I$, $j=1,\ldots,J$, where $M\geq 1$ is a constant. Assume that $G_i=0$ and $H_j=0$ in $B_{\delta_0}(p)\cap O_-$, $i=1,\ldots,I$, $j=1,\ldots,J$. Assume also that $f$ is strongly pseudo-convex at $p$, in the sense of Definition \ref{psconvexqual}, and $L(f)(p)\neq 0$. Then $G_i=0$ and $H_j=0$ in $B_{\delta_1}(p)$, $i=1,\ldots,I$, $j=1,\ldots,J$, for some constant $\delta_1\in(0,\delta_0)$ sufficiently small.
\end{lemma}

Lemma \ref{extendedCarl2} is proved in \cite[Lemma 3.4]{AlIoKl}, using two Carleman inequalities: Proposition 3.3 in \cite{IoKl} and Lemma A.3 in \cite{AlIoKl}. The implicit constant $\delta_1>0$ depends only on constants $A$ in \eqref{quant2}, $\delta_0$, and the constant $A_1$ in the following quantitative form of strong pseudo-convexity:

\begin{lemma}\label{psconvexquan}
(a) Assume that $f$ is strongly pseudo-convex at $p$. Then there are constants $A_1\geq A$ and $\mu\in[-A_1,A_1]$ such that, for any vector $X=X^\al\partial_\al$ at $p$,
\begin{equation}\label{quant6}
\begin{split}
&|\partial^1f(p)|\geq A_1^{-1},\\
&X^\al X^\be(\mu\g_{\al\be}(p)-\D_\al\D_\be f(p))+A_1|X(f)(p)|^2\geq A_1^{-1}|X|^2,
\end{split}
\end{equation}
where $|X|^2=(X^1)^2+\ldots+(X^d)^2$.

(b) Moreover, the inequalities \eqref{quant6} persist in a small neighborhood of $p$, in the sense that there is $\eps_1=\eps(A_1)>0$ such that for any vector-field $X=X^\al\partial_\al$ in $B_{\eps_1}(p)$, the inequalities
\begin{equation}\label{quant7}
\begin{split}
&|\partial^1f|\geq (2A_1)^{-1},\\
&X^\al X^\be(\mu\g_{\al\be}-\D_\al\D_\be f)+A_1|X(f)|^2\geq (2A_1)^{-1}|X|^2,
\end{split}
\end{equation}
hold in $B_{\eps_1}(p)$, where $|X|^2=(X^1)^2+\ldots+(X^d)^2$ and $\mu$ is as in \eqref{quant6}.
\end{lemma}

\begin{proof}[Proof of Lemma \ref{quant6}] (a) The first inequality in \eqref{quant6} is just a quantitative form of the assumption that $p$ is not a critical point of $f$. To derive the second inequality, let $h_{\al\be}=-\D_\al\D_\be f(p)$ and
\begin{equation*}
\delta_0=\inf_{|X|=1,X^\al X_\al=X^\al\D_\al f=0}X^\al X^\be h_{\al\be}.
\end{equation*}
By compactness, this infimum is attained, and it follows from \eqref{qual1} that $\delta>0$. By homogeneity, it follows that
\begin{equation}\label{qual2}
X^\al X^\be h_{\al\be}\geq\delta_0 |X|^2\quad\text{ if }X^\al X_\al=X^\al\D_\al f=0.
\end{equation}

We would like to prove now that there is $n_0\in\{1,2,\ldots\}$ such that
\begin{equation}\label{qual3}
X^\al X^\be h_{\al\be}+n_0(X^\al\D_\al f(p))^2\geq(\delta_0/2) |X|^2\quad\text{ if }X^\al X_\al=0.
\end{equation}
Indeed, otherwise for any $n=1,2,\ldots$ there would exist a vector $X_n=X_n^\al\partial_\al$ such that $|X_n|=1$, $\g_p(X_n,X_n)=0$, and
\begin{equation*}
X_n^\al X_n^\be h_{\al\be}+n(X_n^\al\D_\al f(p))^2\leq\delta_0/2.
\end{equation*}
After passing to a subsequence, we may assume that $X_n$ converges to a vector $X$, with $|X|^2=1$, $X^\al X_\al=0$, $X^\al\D_\al f(p)=0$, and $X^\al X^\be h_{\al\be}\leq\delta_0/2$, which contradicts \eqref{qual2}. Therefore \eqref{qual3} holds for some constant $n_0$.

Let $\mathcal{C}_+=\{X\in T_p\M:|X|=1\text{ and }X^\al X_\al>0\}$, $\mathcal{C}_-=\{X\in T_p\M:|X|=1\text{ and }X^\al X_\al<0\}$, and, for $\delta\in [0,1]$, $\mathcal{C}_\delta=\{X\in T_p\M:|X|=1\text{ and }|X^\al X_\al|\leq\delta\}$. Since the metric $\g$ is non-degenerate, we may assume that $\mathcal{C}_+\neq\emptyset$ (if $\mathcal{C}_+=\emptyset$ then $\mathcal{C}_-\neq\emptyset$ and the proof proceeds in a similar way). For $\rho\in\mathbb{R}$, we consider the function
\begin{equation*}
K_\rho:T_p\M\to\mathbb{R},\qquad K_\rho(X)=X^\al X^\be h_{\al\be}+n_0(X^\al\D_\al f(p))^2+\rho X^\al X_\al,
\end{equation*}
where $n_0$ is as in \eqref{qual3}. Using a simple compactness argument as before, it follows from \eqref{qual3}, that
\begin{equation}\label{hj1}
\text{ there is }\delta'>0\text{ such that }K_0(X)\geq\delta_0/4\text{ for any }X\in\mathcal{C}_{\delta'}.
\end{equation}
Then it follows that there is $\rho_1\geq  0$ sufficiently large such that
\begin{equation*}
K_{\rho_1}\geq 0\text{ if }X\in\mathcal{C}_+\text{ and there is }X\in\mathcal{C}_+\text{ such that }K_{-\rho_1}(X)<0.
\end{equation*}

Let
\begin{equation*}
\rho_0=\inf\{\rho\in[-\rho_1,\rho_1]:K_\rho(X)\geq 0\text{ for any }X\in\mathcal{C}_+\}.
\end{equation*}
We analyze the function $K_{\rho_0}(X)=X^\al X^\be k_{\al\be}$, where
\begin{equation*}
k_{\al\be}:=h_{\al\be}+n_0\D_\al f(p)\D_\be f(p)+\rho_0\g_{\al\be}.
\end{equation*}
In view of the definition of $\rho_0$, $K_{\rho_0}(X)\geq 0$ in $\mathcal{C}_+$. Moreover, using also \eqref{hj1}, there is $X_0\in \mathcal{C}_+$ such that $K_{\rho_0}(X_0)=0$. Since $K_{\rho_0}$ is homogeneous of degree $2$, it follows that the point $X_0$ is a local minimum for $K_{\rho_0}$ in $T_p\M$. Therefore
\begin{equation}\label{hj2}
V^\al X_0^\be k_{\al\be}=0\text{ and }V^\al V^\be k_{\al\be}\geq 0\text{ for any }V\in T_p\OO.
\end{equation}

We show now that
\begin{equation}\label{hj3}
K_{\rho_0}(X)\neq 0\text{ for any }X\in\mathcal{C}_-.
\end{equation}
Indeed, assuming $K_{\rho_0}(X_1)=0$ for some $X_1\in\mathcal{C}_-$, it follows from \eqref{hj2} that $K_{\rho_0}(tX_0+(1-t)X_1)=0$ for any $t\in[0,1]$. However, this contradicts \eqref{hj1} since there is $t_0\in[0,1]$ such that $\g_p(t_0X_0+(1-t_0)X_1,t_0X_0+(1-t_0)X_1)=0$ and $t_0X_0+(1-t_0)X_1\neq 0$.

Using \eqref{hj1}, \eqref{hj2}, and \eqref{hj3} it follows that $K_{\rho_0}(X)>0$ for any $x\in \mathcal{C}_-\cup\mathcal{C}_{\delta''}$, for some $\delta''>0$. A simple compactness argument then shows that there is $n_1$ large enough such that $K_{\rho_0+1/n_1}>0$ in $\{X\in T_p\M:|X|=1\}$. The second inequality in \eqref{quant6} follows by setting $\mu=\rho_0+1/n_1$ and $A_1$ sufficiently large.

Part (b) of the lemma follows from part (a) and the assumption \eqref{quant2}.
\end{proof}

\section{Proof of Theorem \ref{mkathm}}\label{mka}

The plan of the proof is the following: we fix a point $p\in U_0\cap\mathcal{H}^-\cap\overline{\mathbf{E}}$, outside both the bifurcation sphere $S_0=\mathcal{H}^-\cap\mathcal{H}^+$ and the axis of symmetry $\mathcal{A}=\{p\in\overline{\mathbf{E}}:Z(p)=0\}$. Then we consider the Kerr metric $\g$ and the induced metric
\begin{equation*}
h_{\al\be}=X\g_{\al\be}-\T_\al \T_\be,\qquad\text{ where }X=\g(\T,\T),
\end{equation*}
on a hypersurface $\Pi$ passing through the point $p$ and transversal to $\T$. The metric $h$ is nondegenerate (Lorentzian) as long as $X>0$ in $\Pi$, which explains our assumption $0<a<m$. It is well-known, see for example \cite[Section 3]{We}, that the Einstein vacuum equations together with stationarity $\mathcal{L}_\T\g=0$ are equivalent to the system of equations
\begin{equation}\label{brac}
\begin{split}
&{}^h\mathbf{Ric}_{ab}=\frac{1}{2X^2}(\nabla_a X\nabla_b X+\nabla_a Y\nabla_b Y),\\
&{}^h\square (X+iY)=\frac{1}{X}h^{ab}\partial_a(X+iY)\partial_b(X+iY),
\end{split}
\end{equation}
in $\Pi$, where $Y$ is the Ernst potential associated to $\T$. We rederive these equations in Proposition \ref{kerex} below, together with other explicit equations and identities that are needed for the proof of the theorem.

We then modify the metric $h$ and the functions $X$ and $Y$ in a neighborhood of the point $p$ in such a way that the identities \eqref{brac} are still satisfied. The existence of a large family of smooth triplets $(\widetilde{h},\widetilde{X},\widetilde{Y})$ satisfying \eqref{brac} and agreeing with the Kerr triplet in $\Pi\setminus\mathbf{E}$ follows by solving a characteristic initial-value problem, using the main theorem in \cite{Re}. 

Finally, in Proposition \ref{mainprop2} we construct the space-time metric $\widetilde{\g}$,
\begin{equation*}
\widetilde{\g}_{ab}=\widetilde{X}^{-1}\widetilde{h}_{ab}+\widetilde{X}\widetilde{A}_a\widetilde{A}_b,\qquad \widetilde{\g}_{a4}=\widetilde{X}\widetilde{A}_a,\qquad\widetilde{\g}_{44}=\widetilde{X},\qquad a,b=1,2,3,
\end{equation*}
 associated to the triplet $(\widetilde{h},\widetilde{X},\widetilde{Y})$, the vector-field $\T=\partial_4$, and a suitable $1$-form $\widetilde{A}$ which is defined in $\Pi$. By construction and \cite[Theorem 1]{We}, this metric verifies the identities ${}^{\widetilde{\g}}\mathbf{Ric}=0$ and $\Lie_\T\widetilde{\g}=0$, in a suitable open set $U$. Then we show that we have enough flexibility to choose initial conditions for $\widetilde{X},\widetilde{Y}$ such that the vector-field $Z$ cannot be extended as a Killing vector-field for $\widetilde{\g}$ commuting with $\T$, in the open set $U$.

\subsection{Explicit calculations}

We consider the Kerr space-time $\mathcal{K}(m,a)$ in standard Boyer--Lindquist coordinates,
\begin{equation}\label{k1}
\g=-\frac{q^2\Delta}{\Sigma^2}(dt)^2+\frac{\Sigma^2(\sin\theta)^2}{q^2}\Big(d\phi-\frac{2amr}{\Sigma^2}dt\Big)^2 +\frac{q^2}{\Delta}(dr)^2+q^2(d\theta)^2,
\end{equation}
where
\begin{equation}\label{k2}
\begin{cases}
&\Delta=r^2+a^2-2mr;\\
&q^2=r^2+a^2(\cos\theta)^2;\\
&\Sigma^2=(r^2+a^2)q^2+2mra^2(\sin\theta)^2=(r^2+a^2)^2-a^2(\sin\theta)^2\Delta.
\end{cases}
\end{equation}

We make the change of variables
\begin{equation*}
du_-=dt-(r^2+a^2)\Delta^{-1}dr,\qquad d\phi_-=d\phi-a\Delta^{-1}dr.
\end{equation*}
In the new coordinates $(\theta,r,\phi_-,u_-)$ the space-time metric becomes
\begin{equation}\label{gro1}
\begin{split}
\g&=q^2d\theta^2-(du_-\otimes dr+dr\otimes du_-)+a(\sin\theta)^2(d\phi_-\otimes dr+dr\otimes d\phi_-)\\
&-\frac{2amr(\sin\theta)^2}{q^2}(d\phi_-\otimes du_-+du_-\otimes d\phi_-)+\frac{\Sigma^2(\sin\theta)^2}{q^2}d\phi_-^2+\frac{2mr-q^2}{q^2}du_-^2,
\end{split}
\end{equation}
and the vector-field $\T=d/dt$ becomes $\T=d/du_-$. The metric $\g$ and the vector-field $\T$ are smooth in the region
\begin{equation*}
R=\{(\theta,r,\phi_-,u_-)\in(0,\pi)\times(0,\infty)\times(-\pi,\pi)\times\mathbb{R}:2mr-q^2>0\}.
\end{equation*}

Let
\begin{equation*}
X=\g(\T,\T)=\frac{2mr-q^2}{q^2},\qquad h_{\al\be}=X\g_{\al\be}-\T_\al\T_\be,
\end{equation*}
and
\begin{equation*}
\Pi=\{(\theta,r,\phi_-,u_-)\in R:u_-=0\}.
\end{equation*}
Let
\begin{equation}\label{vectors}
\partial_1=\frac{d}{d\theta},\qquad\partial_2=\frac{d}{dr},\qquad\partial_3=\frac{d}{d\phi_-},
\end{equation}
denote the vector-fields in $\Pi$ induced by coordinates $(\theta,r,\phi_-)$. We calculate the components of the metric $h$ along the surface $\Pi$,
\begin{equation}\label{no1}
\begin{split}
&h_{11}=2mr-q^2,\qquad h_{12}=0,\qquad h_{13}=0,\qquad h_{22}=-1,\\
&h_{23}=-a(\sin\theta)^2,\qquad h_{33}=-\Delta(\sin\theta)^2.
\end{split}
\end{equation}
Therefore
\begin{equation}\label{no2}
\begin{split}
&h^{11}=\frac{1}{2mr-q^2},\qquad h^{12}=0,\qquad h^{13}=0,\qquad h^{22}=\frac{\Delta}{2mr-q^2},\\
&h^{23}=\frac{-a}{2mr-q^2},\qquad h^{33}=\frac{1}{(\sin\theta)^2(2mr-q^2)}.
\end{split}
\end{equation}

Let
\begin{equation}\label{gro6}
\Gamma_{cab}=h(\nabla_{\partial_b}\partial_a,\partial_c)=(1/2)(\partial_a h_{bc} +\partial_b h_{ac}-\partial_c h_{ab}),\qquad {\Gamma^d}_{ab}=h^{cd}\Gamma_{cab}.
\end{equation}
Using \eqref{no1} and \eqref{no2} we calculate
\begin{equation}\label{no3}
\begin{split}
&{\Gamma^1}_{11}=\frac{a^2\sin\theta\cos\theta}{2mr-q^2},\qquad {\Gamma^2}_{11}=\frac{\Delta(r-m)}{2mr-q^2},\qquad {\Gamma^3}_{11}=\frac{a(m-r)}{2mr-q^2},\\
&{\Gamma^1}_{12}=\frac{m-r}{2mr-q^2},\qquad {\Gamma^2}_{12}=\frac{a^2\sin\theta\cos\theta}{2mr-q^2},\qquad {\Gamma^3}_{12}=\frac{-a\cot\theta}{2mr-q^2},\\
&{\Gamma^1}_{13}=0,\qquad {\Gamma^2}_{13}=0,\qquad {\Gamma^3}_{13}=\cot\theta,\\
&{\Gamma^1}_{22}=0,\qquad {\Gamma^2}_{22}=0,\qquad {\Gamma^3}_{22}=0,\\
&{\Gamma^1}_{23}=\frac{a\sin\theta\cos\theta}{2mr-q^2},\qquad {\Gamma^2}_{23}=\frac{a(r-m)(\sin\theta)^2}{2mr-q^2},\qquad {\Gamma^3}_{23}=\frac{m-r}{2mr-q^2},\\
&{\Gamma^1}_{33}=\frac{\Delta\sin\theta\cos\theta}{2mr-q^2},\qquad {\Gamma^2}_{33}=\frac{\Delta(r-m)(\sin\theta)^2}{2mr-q^2},\qquad {\Gamma^3}_{33}=\frac{a(m-r)(\sin\theta)^2}{2mr-q^2}.
\end{split}
\end{equation}

We have
\begin{equation*}
\begin{split}
{}^{h}\R(\partial_a,\partial_c)\partial_b&=\nabla_{\partial_a}(\nabla_{\partial_c}\partial_b)- \nabla_{\partial_c}(\nabla_{\partial_a}\partial_b)\\
&=\nabla_{\partial_a}({\Gamma^d}_{bc}\partial_d)-\nabla_{\partial_c}({\Gamma^d}_{ba}\partial_d)\\
&={\partial_a}({\Gamma^d}_{bc})\partial_d+{\Gamma^d}_{bc}{\Gamma^e}_{da}\partial_e-{\partial_c}({\Gamma^d}_{ba})\partial_d +{\Gamma^d}_{ba}{\Gamma^e}_{dc}\partial_e\\
&=[{\partial_a}({\Gamma^d}_{bc})-{\partial_c}({\Gamma^d}_{ba})+{\Gamma^e}_{bc}{\Gamma^d}_{ea} -{\Gamma^e}_{ba}{\Gamma^d}_{ec}]\partial_d,
\end{split}
\end{equation*}
therefore
\begin{equation}\label{gro7}
{}^h\mathbf{Ric}_{ba}={\partial_c}({\Gamma^c}_{ba})-{\partial_a}({\Gamma^c}_{bc})+{\Gamma^d}_{ba}{\Gamma^c}_{dc} -{\Gamma^d}_{cb}{\Gamma^c}_{da}.
\end{equation}
Using \eqref{gro6} we calculate
\begin{equation}\label{gro8}
{\Gamma^c}_{bc}=(1/2)h^{ca}(\partial_{b}h_{ca})=(1/2)\partial_b(\log |h|)=\partial_b(\log(\sin\theta(2mr-q^2))).
\end{equation}
Thus
\begin{equation}\label{Rich}
\begin{split}
&{}^h\mathbf{Ric}_{11}=\frac{2m^2a^2(\sin\theta)^2}{(2mr-q^2)^2},\qquad{}^h\mathbf{Ric}_{12}=0, \qquad{}^h\mathbf{Ric}_{13}=0,\\
&{}^h\mathbf{Ric}_{22}=\frac{2m^2}{(2mr-q^2)^2},\qquad{}^h\mathbf{Ric}_{23}=0,\qquad{}^h\mathbf{Ric}_{33} =0.
\end{split}
\end{equation}

Let
\begin{equation}\label{pote}
\begin{split}
&X=\frac{2mr-q^2}{q^2},\qquad Y=-\frac{2ma\cos\theta}{q^2},\\
&T_{ab}=\frac{1}{2X^2}(\nabla_a X\nabla_b X+\nabla_a Y\nabla_b Y).
\end{split}
\end{equation}
We calculate
\begin{equation}\label{pote2}
\begin{split}
&\partial_1X=\frac{4a^2mr\sin\theta\cos\theta}{q^4},\qquad\partial_2X=\frac{2mq^2-4mr^2}{q^4},\qquad\partial_3X=0,\\
&\partial_1 Y=\frac{2ma\sin\theta q^2-4ma^3\sin\theta(\cos\theta)^2}{q^4},\qquad\partial_2 Y=\frac{4mra\cos\theta}{q^4},\qquad \partial_3Y=0.
\end{split}
\end{equation}
Therefore
\begin{equation*}
\begin{split}
T_{11}=\frac{2m^2a^2(\sin\theta)^2}{(2mr-q^2)^2},\quad T_{12}=0,\quad T_{13}=0,\quad T_{22}=\frac{2m^2}{(2mr-q^2)^2},\quad T_{23}=0,\quad T_{33}=0.
\end{split}
\end{equation*}
Using also \eqref{Rich} it follows that
\begin{equation*}
{}^h\mathbf{Ric}=T.
\end{equation*}

Using \eqref{no2}, \eqref{pote2}, and $|h|=(\sin\theta)^2(2mr-q^2)^2$ we calculate
\begin{equation}\label{lol1}
\begin{split}
&|h|^{1/2}h^{1j}\partial_j(X+iY)=i\frac{2am(\sin\theta)^2}{q^4}(r-ia\cos\theta)^2,\\
&|h|^{1/2}h^{2j}\partial_j(X+iY)=\frac{-2m\Delta\sin\theta}{q^4}(r-ia\cos\theta)^2,\\
&|h|^{1/2}h^{3j}\partial_j(X+iY)=\frac{2ma\sin\theta}{q^4}(r-ia\cos\theta)^2.
\end{split}
\end{equation}
Therefore
\begin{equation}\label{box}
\begin{split}
&{}^h\square X=\frac{24m^2r^2a^2(\cos\theta)^2-4m^2r^4-4m^2a^4(\cos\theta)^4} {q^6(2mr-q^2)},\\
&{}^h\square Y=\frac{16m^2ra\cos\theta(r^2-a^2(\cos\theta)^2)}{q^6(2mr-q^2)}.
\end{split}
\end{equation}
We calculate also
\begin{equation*}
\begin{split}
&X^{-1}h^{ij}(\partial_iX\partial_jX-\partial_iY\partial_jY)=\frac{24m^2r^2a^2(\cos\theta)^2-4m^2r^4-4m^2a^4(\cos\theta)^4} {q^6(2mr-q^2)},\\
&2X^{-1}h^{ij}\partial_iX\partial_jY=\frac{16m^2ra\cos\theta(r^2-a^2(\cos\theta)^2)}{q^6(2mr-q^2)}.
\end{split}
\end{equation*}
Therefore
\begin{equation}\label{system}
\begin{split}
&{}^h\mathbf{Ric}_{ab}=\frac{1}{2X^2}(\nabla_a X\nabla_b X+\nabla_a Y\nabla_b Y),\\
&{}^h\square (X+iY)=\frac{1}{X}h^{ab}\partial_a(X+iY)\partial_b(X+iY).
\end{split}
\end{equation}

The components of the  spacetime  metric $\g$  in the coordinates 
$(\theta,r,\phi_-,u_-)$ (see \eqref{gro1}) have the form,
\begin{equation*}
{\g}_{ab}={X}^{-1}{h}_{ab}+{X}{A}_a{A}_b,\qquad {\g}_{a4}={X}{A}_a,\qquad{\g}_{44}={X},\qquad a,b=1,2,3.
\end{equation*} 
or, with $x=(\th, r, \phi)$,
\begin{equation*}
\g=(X du_-  +A_a dx^a)^2 +X^{-1}h_{ab}dx^a dx^b   \label{gro1-h}
\end{equation*}
where,
\begin{equation}\label{as1}
A_1=0,\qquad A_2=-\frac{q^2}{2mr-q^2},\qquad A_3=-\frac{2amr(\sin\theta)^2}{2mr-q^2},
\end{equation}
We compute
\begin{equation*}
\begin{split}
&\partial_1A_2-\partial_2A_1=\frac{4a^2mr\sin\theta\cos\theta}{(2mr-q^2)^2},\\
&\partial_2A_3-\partial_3A_2=\frac{-2ma(\sin\theta)^2 (r^2-a^2(\cos\theta)^2)}{(2mr-q^2)^2},\\
&\partial_3A_1-\partial_1A_3=\frac{-4mra\Delta\sin\theta\cos\theta}{(2mr-q^2)^2}.
\end{split}
\end{equation*}
Therefore, using also \eqref{lol1}, with ${}^h\negmedspace\in_{123}=-|h|^{1/2}$,
\begin{equation}\label{as2}
X^2(\nabla_aA_b-\nabla_bA_a)={}^h\negmedspace\in_{abc}\nabla^cY.
\end{equation}

To summarize, we verified the following:

\begin{proposition}\label{kerex}
With the notation above, the metric $h$, the functions $X,Y$, and the $1$-form $A$ satisfy the identities (in $\Pi$)
\begin{equation*}
\begin{split}
&{}^h\mathbf{Ric}_{ab}=\frac{1}{2X^2}(\nabla_a X\nabla_b X+\nabla_a Y\nabla_b Y),\\
&{}^h\square (X+iY)=\frac{1}{X}h^{ab}\partial_a(X+iY)\partial_b(X+iY),\\
&X^2(\nabla_aA_b-\nabla_bA_a)={}^h\negmedspace\in_{abc}\nabla^cY.
\end{split}
\end{equation*}
\end{proposition}
\begin{remark} Under a change of coordinates of the form $u'_-=u_--f(x^1,x^2,x^3)$
the $1$-form $A=A_a dx^a$ changes according to the formula $A'=A-d f$.
 The change of coordinates amounts to a  choice  of the  hypersurface $\Pi$ i.e. 
 instead of $u_-=0$ we would chose $u_-=f(\th, r, \phi)$.
  \label{rem:changeA}
\end{remark}

\subsection{The metric $\widetilde{h}$}\label{mem6} We would like to construct now a large family of triplets $(\widetilde{h},\widetilde{X},\widetilde{Y})$ and $1$-forms $\widetilde{A}$, such that the identities in Proposition \ref{kerex} are still satisfied in a neighborhood in $\Pi$ of a fixed point $p\in (U_0\cap\mathcal{H}^-\cap\overline{\mathbf{E}})\setminus(\mathcal{A}\cup S_0)$. Let
\begin{equation*}
\mathcal{N}_0=\{x\in S:r(x)=r_+:=m+\sqrt{m^2-a^2}\}.
\end{equation*}
This is a 2-dimensional hypersurface in $\Pi$; the vector-fields $\partial_1$ and $\partial_3$ are tangent to $\mathcal{N}_0$ and, using \eqref{no1} and \eqref{no3},
\begin{equation*}
h(\partial_3,\partial_3)=h(\partial_3,\partial_1)=0,\quad \nabla_{\partial_3}\partial_3=-[(m/a)^2-1]^{1/2}\partial_3,\qquad \text{ along } \mathcal{N}_0.
\end{equation*}
Therefore $\mathcal{N}_0$ is a null hypersurface in $\Pi$. Along $\mathcal{N}_0\subset \Pi$  we define the smooth, transversal, null  vector-field,
\begin{equation}\label{L1}
L=(2a^2(\sin\theta)^2-\Delta)^{-1}\cdot [2a\partial_2-(\sin\theta)^{-2}\partial_3].
\end{equation}
Using \eqref{no1}, it follows that
\begin{equation}\label{L2}
h(L,L)=h(L,\partial_1)=0,\qquad [L,\partial_3]=0,\qquad h(L,\partial_3)=-1,\qquad \text{ along }\mathcal{N}_0.
\end{equation}

Let
\begin{equation*}
P=\{x\in\mathcal{N}_0:\phi_-(x)=0\},\qquad p=\{x\in P:\theta(x)=\theta_0\in(0,\pi)\}.
\end{equation*}
Thus $P$ is a 1-dimensional smooth curve in $\mathcal{N}_0$ and $p\in P$ is a point. We extend the vector-field $L$ to a small open neighborhood $D$ of $p$ in $\Pi$, by solving the geodesic equation
\begin{equation*}
\nabla_LL=0 \qquad\text{ in }D.
\end{equation*}
Then we construct the null hypersurface $\mathcal{N}_1$ in $D$ as the congruence of geodesic curves tangent to $L$ and passing through the curve $P$. We also fix a time-orientation in $D$ such that $\partial_3$ and $L$ are future-directed null vector-fields along $P\cap D$. and let $J^+(\mathcal{N}_1)$ denote the causal future of $\mathcal{N}_1$ in $D$. Let
\begin{equation*}
D_-=\{x\in D:\Delta(x)<0\},\qquad D_+=\{x\in D:\Delta(x)>0\}.
\end{equation*}
The following proposition is a consequence of the main theorem in \cite{Re}.

\begin{proposition}\label{mainprop}
Assume $\widetilde{X},\widetilde{Y}:\mathcal{N}_1\to\mathbb{R}$ are smooth functions satisfying
\begin{equation*}
\widetilde{X}=X\text{ and }\widetilde{Y}=Y\text{ in }\mathcal{N}_1\cap D_-.
\end{equation*}
Then there is a small neighborhood $D'$ of $p$ in $\Pi$, a smooth metric $\widetilde{h}$ in $J^+(\mathcal{N}_1)\cap D'$,  and smooth extensions $\widetilde{X},\widetilde{Y}:J^+(\mathcal{N}_1)\cap D'\to\mathbb{R}$ such that, in $J^+(\mathcal{N}_1)\cap D'$,
\begin{equation}\label{cond1}
\begin{split}
&{}^{\widetilde{h}}\mathbf{Ric}_{ab}=\frac{1}{2\widetilde{X}^2}(\nabla_a \widetilde{X}\nabla_b \widetilde{X}+\nabla_a \widetilde{Y}\nabla_b \widetilde{Y}),\\
&{}^{\widetilde{h}}\square (\widetilde{X}+i\widetilde{Y})=\frac{1}{\widetilde{X}}\widetilde{h}^{ab}\partial_a(\widetilde{X}+i\widetilde{Y}) \partial_b(\widetilde{X}+i\widetilde{Y}).
\end{split}
\end{equation}
In addition
\begin{equation}\label{cond2}
\widetilde{X}=X,\quad\widetilde{Y}=Y,\quad\widetilde{h}=h\qquad\text{ in }
J^+(\mathcal{N}_1)\cap D'\cap D_-,
\end{equation}
and, for any vector-field $V$ tangent to $\mathcal{N}_1\cap D'$,
\begin{equation}\label{cond3}
\widetilde{h}(L,V)=0\text{ and }\widetilde{\nabla}_LL=0\qquad\text{ along }\mathcal{N}_1\cap D'.
\end{equation}
\end{proposition}

To be able to construct the desired space-time metric $\widetilde{\g}$ we also need to extend the $1$-form $A$ (compare with the formula \eqref{spa1}). More precisely:

\begin{proposition}\label{Aextension}
There is a smooth $1$-form $\widetilde{A}_a$ in a neighborhood $\mathcal{D}$ of $p$ in $J^+(\mathcal{N}_1)$ satisfying (compare with \eqref{as2})
\begin{equation}\label{as3}
\begin{split}
&\widetilde{X}^2(\widetilde{\nabla}_a\widetilde{A}_b-\widetilde{\nabla}_b\widetilde{A}_a)=\widetilde{\in}_{abc}\widetilde{\nabla}^c\widetilde{Y},\\
&\widetilde{A}=A\qquad\text{ in }\mathcal{D}\cap D_-.
\end{split}
\end{equation}
\end{proposition}
\begin{proof}
Without loss of generality\footnote{Alternatively the argument below can easily be adapted 
 to the case $L\c A\neq0 $ by a straightforward modification of  equation \eqref{as4}.    } we may assume that  $L^a A_a=0$ in  a full neighborhood $D$ of $p$ in $\Pi$.
Indeed in  view of remark \ref{rem:changeA}     we can choose  a function $f$ in $D$ such that $L(f)=L^a A_a$ and  change  $\Pi$ to $\Pi'$ by  redefining   $u_-'=u_--f$.  In $\Pi'$ the corresponding $L', A'$ 
verify $(L')^a\c A'_a=0$.

Let $\widetilde{L}$ denote the geodesic vector-field (i.e. $\widetilde{\nabla}_{\widetilde{L}}\widetilde{L}=0$) generated in a neighborhood of the point $p$ in by the vector-field $L$ define on $\mathcal{N}_0$ in \eqref{L1}, so
\begin{equation*}
\widetilde{L}=L\text{ in }[(J^+(\mathcal{N}_1)\cap D_-)\cup\mathcal{N}_1]\cap D'.
\end{equation*}

We then define the form $\widetilde{A}$ as the solution of the transport equation, in a neighborhood of the point $p$ in $J^+(\mathcal{N}_1)$,
\begin{equation}\label{as4}
\begin{split}
&\widetilde{L}^a\widetilde{\nabla}_a\widetilde{A}_b+\widetilde{A}_a\widetilde{\nabla}_b\widetilde{L}^a= \widetilde{\in}_{abc}\widetilde{X}^{-2}\widetilde{L}^a\widetilde{\nabla}^c\widetilde{Y},\\
&\widetilde{A}=A\qquad\text{ along }\mathcal{N}_0,
\end{split}
\end{equation}
 It follows easily from \eqref{as2} that the form $A$ verifies this transport equation in $D_-$, thus $\widetilde{A}$ is a well-defined smooth form in a neighborhood $\mathcal{D}$ of $p$ in $J^+(\mathcal{N}_1)$ and $\widetilde{A}=A$ in $D_-$. It remains to prove the identity in the first line of \eqref{as3}. We observe first
 that $\widetilde {A}_a \widetilde{L}^a$   vanishes in a neighborhood of  $p$ in    $J^+(\mathcal{N}_1)$. Indeed,
 using the definition \eqref{as4} we compute,
\begin{equation*}
\widetilde{L}^a\widetilde{\nabla}_a(\widetilde{L}^b\widetilde{A}_b) =\widetilde{L}^b\widetilde{L}^a\widetilde{\nabla}_a\widetilde{A}_b=0.
\end{equation*}
therefore
\begin{equation}\label{as5}
\widetilde{L}^b\widetilde{A}_b= 0\qquad\text{ in a neighborhood of }p\text{ in }J^+(\mathcal{N}_1).
\end{equation}
Letting
\begin{equation}\label{as5.5}
\widetilde{Q}_{ab}=\widetilde{X}^2(\widetilde{\nabla}_a\widetilde{A}_b-\widetilde{\nabla}_b\widetilde{A}_a) -\widetilde{\in}_{abc}\widetilde{\nabla}^c\widetilde{Y},
\end{equation}
it follows from \eqref{as4} and \eqref{as5} that
\begin{equation}\label{as7}
\widetilde{L}^a\widetilde{Q}_{ab}=0,\qquad \widetilde{L}^b\widetilde{Q}_{ab}=0.
\end{equation}
To show that $\widetilde{Q}$ vanish identically we derive a transport equation for it. In fact we 
 show   in the lemma below  that  $\Lie_{\widetilde{L}} (\widetilde{X}^{-2} \widetilde{Q} )$ vanishes identically   in a neighborhood of $p$ in
 $J^+(\mathcal{N}_1)$.  Since $\widetilde{Q}$ vanishes in $D_-$  it follows that $\widetilde{Q}$
  vanishes in a neighborhood of $p$ in $J^+(\mathcal{N}_1)$, as desired. Thus the proof reduces to the lemma below.
  \end{proof}

\begin{lemma}\label{seri}
Consider a $3$-dimensional  Lorentzian manifold $(\Pi, h)$   and scalar functions ($ X , Y)$ which  verify  the equation,
\bea
\square _h Y=2 X^{-1}h^{ab}\nabla_a X\nabla _b Y,\label{eqY}
\eea
 Assume also given a  $1$-form $A$ which verifies, 
\bea
L^a\nabla_a A_b +A_a \nabla_b L^a=X^{-2}\in_{abc} L^a \nabla^c Y, \qquad  L\c A=0, \label{def-extA}
\eea
with    $L$  a  null geodesic vector-field in $\Pi$. Then the $2$-form
  \bea
Q_{ab}=X^2(\nabla_a A_b-\nabla_b A_a)-\in_{abc}\nabla^c Y
\eea
verifies the equation,
\bea
\Lie_L  (X^{-2}Q)=0.
\eea
\end{lemma} 

\begin{proof}[Proof of Lemma \ref{seri}]
We have,
\beaa
\Lie_L Q_{ab}&=& 2 X\Lie_L X(\nabla_a A_b-\nabla_b A_a)     +X^2   \Lie_L (\nabla_a A_b-\nabla_b A_a)\\
 &-&( \Lie_L \in_{abc})    \nabla^c Y-\in_{abc}\Lie_L (\nabla^c Y)\\
 &=&2 X^{-1}\nabla_L X Q_{ab} +2 X^{-1} \Lie_L X ( \in_{abc}  \nabla^c Y) +X^2  (\nabla_a  \Lie_L A_b-\nabla_b  \Lie_L A_a)    \\
 &-&( \Lie_L \in_{abc})    \nabla^c Y-\in_{abc}\Lie_L (\nabla^c Y)\\
 &=&2 X^{-1}\nabla_L X Q_{ab}+X^2  (\nabla_a  \Lie_L A_b-\nabla_b  \Lie_L A_a) 
 +\in_{abc}\big( 2X^{-1}\nabla_L X -\div(L) )\nabla^c Y\\
 &-&\in_{abc}(L^d\nabla_d\nabla^c Y-\nabla^d Y \nabla_d L^c) 
 \eeaa
 Using the equation for $A$, written in the form $\Lie_L A_a=-X^{-2} \in_{acd} L^c \nabla^d Y$,
 \beaa
 X^2  (\nabla_a  \Lie_L A_b-\nabla_b  \Lie_L A_a)&=&
 X^2  (\nabla_a  \Lie_L A_b-\nabla_b  \Lie_L A_a)\\
 &=&
2X^{-1}\big( \nabla_a X\in_{bcd}-\nabla_b X\in_{acd}     \big) L^c \nabla^d Y\\
 &-&\in_{bcd}\nabla_a L^c\nabla^d Y     +\in_{acd} \nabla_bL^c \nabla ^d Y\\
 &-&\in_{bcd}L^c\nabla_a \nabla^d Y+\in_{acd}L^c\nabla_b \nabla^d Y.
 \eeaa 
 Hence,
  \beaa
 \Lie_L Q_{ab}-2 X^{-1}\nabla_L X Q_{ab}&=&E_{ab}
 \eeaa
 with
\beaa
 E_{ab}&=&-\in_{bcd}L^c\nabla_a \nabla^d Y+\in_{acd}L^c\nabla_b \nabla^d Y-\in_{abc}L^d\nabla_d\nabla^c Y\\
 &+&2X^{-1}\big( \nabla_a X\in_{bcd}-\nabla_b X\in_{acd}     \big) L^c \nabla^d Y+
  \in_{abc}2X^{-1}\nabla_L X\nabla^c Y\\
   &-&\in_{bcd}\nabla_a L^c\nabla^d Y     +\in_{acd} \nabla_bL^c \nabla ^d Y 
  +  \in_{abc}\big(  -(\div L) \nabla^c Y+\nabla^d Y \nabla_d L^c\big).
 \eeaa
 To  check that $E\equiv 0$ it suffices to show that its Hodge  dual $\dual E_m:=\frac 1 2 \in_{m}\,^{ab} E_{ab}$
  vanishes.   By a straightforward calculation, involving the usual rules of contracting   tensor products of    the volume form 
      $\in$, we find,
  \beaa
  \dual E_m&=& (\square_h Y-2X^{-1} \nabla_a X\nabla^a Y)L_m \\  
  \eeaa
  from which the lemma  easily follows.
 \end{proof}

\subsection{The space-time metric} Let $\widetilde{X}$, $\widetilde{Y}$, $\widetilde{h}$, $\mathcal{D}$, and $\widetilde{A}$ be as before. In $\mathcal{D}\times I$, where $I\subset \mathbb{R}$ is an open interval, we define the Lorentz metric $\widetilde{\g}$ by
\begin{equation}\label{spa1}
\widetilde{\g}_{ab}=\widetilde{X}^{-1}\widetilde{h}_{ab}+\widetilde{X}\widetilde{A}_a\widetilde{A}_b,\qquad \widetilde{\g}_{a4}=\widetilde{X}\widetilde{A}_a,\qquad\widetilde{\g}_{44}=\widetilde{X},\qquad a,b=1,2,3.
\end{equation}
The functions $\widetilde{X},\widetilde{Y},\widetilde{A}_a,\widetilde{h}_{ab}$, originally defined in $\mathcal{D}$ are extended to $\mathcal{D}\times I$ by
\begin{equation}\label{spa3}
\partial_4(\widetilde{X})=\partial_4(\widetilde{Y})=\partial_4(\widetilde{A}_a)=\partial_4(\widetilde{h}_{ab})=0,\qquad a,b=1,2,3.
\end{equation}
Using \eqref{spa1}, it follows that, with $\widetilde{A}^a=\widetilde{h}^{ab}\widetilde{A}_b$, $a=1,2,3$,
\begin{equation}\label{spa4}
\widetilde{\g}^{ab}=\widetilde{X}\widetilde{h}^{ab},\qquad \widetilde{\g}^{a4}=-\widetilde{X}\widetilde{A}^a,\qquad\widetilde{\g}^{44}=\widetilde{X}^{-1}+\widetilde{X}\widetilde{A}^a\widetilde{A}_a, \qquad|\widetilde{\g}|=\widetilde{X}^{-2}|\widetilde{h}|.
\end{equation}

\begin{proposition}\label{mainprop2}
(a) The metric $\widetilde{\g}$ agrees with the Kerr metric $\g$ in $(\mathcal{D}\cap D_-)\times I$ and satisfies
\begin{equation*}
\Lie_{\partial_4}\widetilde{\g}=0,\qquad {}^{\widetilde{\g}}\mathbf{Ric}=0\qquad\text{ in }\mathcal{D}\times\mathbb{R}.
\end{equation*}
(b) If $Z=Z^4\partial_4+Z^a\partial_a$ is a Killing vector-field for $\widetilde{\g}$ in $\mathcal{D}\times I$ and if $[Z,\partial_4]=0$ then $Z'=Z^a\partial_a$ is a Killing vector-field for $\widetilde{h}$ in $\mathcal{D}$ satisfying $Z'(\widetilde{X})=Z'(\widetilde{Y})=0$, i.e.
\begin{equation}\label{mla0}
Z'(\widetilde{X})=Z'(\widetilde{Y})=0,\qquad(\Lie_{Z'}\widetilde{h})_{ab}=0.
\end{equation}
\end{proposition}

\begin{proof}[Proof of Proposition \ref{mainprop2}] (a) The claims follow easily from definitions, except for
\begin{equation*}
{}^{\widetilde{\g}}\mathbf{Ric}=0\qquad\text{ in }\mathcal{D}\times\mathbb{R}.
\end{equation*}
On the other hand, this is a well-known consequence of the identities \eqref{cond1} and \eqref{as3} satisfied by $\widetilde{h},\widetilde{X},\widetilde{Y}$ and $\widetilde{A}$, and the definitions \eqref{spa1} and \eqref{spa3}. See, for example, \cite[Section 3]{We} for the proof.

(b) The identities $\partial_4Z^4=0$, $\partial_4Z^a=0$, $(\Lie_Z\widetilde{\g})_{44}=0$, $(\Lie_Z\widetilde{\g})_{a4}=0$, and $(\Lie_Z\widetilde{\g})_{ab}=0$ give
\begin{equation*}
\begin{split}
&Z(\widetilde{X})=0,\qquad Z(\widetilde{X}\widetilde{A}_a)+\partial_aZ^c\widetilde{X}\widetilde{A}_c+\partial_aZ^4\widetilde{X}=0,\\
&Z(\widetilde{X}^{-1}\widetilde{h}_{ab}+\widetilde{X}\widetilde{A}_a\widetilde{A}_b)+\partial_aZ^c (\widetilde{X}^{-1}\widetilde{h}_{cb}+\widetilde{X}\widetilde{A}_c\widetilde{A}_b)+\partial_aZ^4\widetilde{X}\widetilde{A}_b\\ &+\partial_bZ^c (\widetilde{X}^{-1}\widetilde{h}_{ac}+\widetilde{X}\widetilde{A}_a\widetilde{A}_c)+\partial_bZ^4\widetilde{X}\widetilde{A}_a=0.
\end{split}
\end{equation*}
Using also \eqref{spa3}, it follows that
\begin{equation*}
\begin{split}
&Z'(\widetilde{X})=0,\qquad Z'(\widetilde{A}_a)+\partial_aZ^c\widetilde{A}_c+\partial_aZ^4=0,\\
&Z'(\widetilde{h}_{ab})+\partial_aZ^c\widetilde{h}_{cb}+\partial_bZ^c\widetilde{h}_{ac}=0.
\end{split}
\end{equation*}
Therefore, along $\mathcal{D}$
\begin{equation*}
Z'(\widetilde{X})=0,\qquad (\Lie_{Z'}\widetilde{h})_{ab}=0,\qquad (\Lie_{Z'}\widetilde{A})_a=-\partial_aZ^4.
\end{equation*}
The last identity in \eqref{mla0}, $Z'(\widetilde{Y})=0$, follows from \eqref{as3}, rewritten in the form
\begin{equation*}
\widetilde{\nabla}^m\widetilde{Y}=-\widetilde{X}^2\widetilde{\in}^{abm}\widetilde{\nabla}_a\widetilde{A}_b.
\end{equation*}
\end{proof}

We can now complete the proof of the theorem.

\begin{proof}[Proof of Theorem \ref{mkathm}] We fix a point $p\in (U_0\cap\mathcal{H}^-\cap\overline{\mathbf{E}})\setminus(\mathcal{A}\cup S_0)$; we may assume that
\begin{equation*}
u_-(p)=0,\qquad \phi_-(p)=0,\qquad \theta(p)\in(0,\pi),\qquad r(p)=m+\sqrt{m^2-a^2}.
\end{equation*}
Then we define the surface $\mathcal{N}_1$ as in Proposition \ref{mainprop}. For any smooth functions $\widetilde{X},\widetilde{Y}:\mathcal{N}_1\to\mathbb{R}$ agreeing with $X,Y$ in $\mathcal{N}_1\cap D_-$, we construct the corresponding neighborhood $\mathcal{D}$ of $p$ in $J^+(\mathcal{N}_1)$ (which we may assume to be diffeomorphic to the unit ball in $\mathbb{R}^3$ and sufficiently small relative to $U_0$), the smooth Lorentzian metric $\widetilde{h}$ in $\mathcal{D}$, the scalars $\widetilde{X},\widetilde{Y}:\mathcal{D}\to\mathbb{R}$, and the $1$-form $\widetilde{A}$, verifying (see \eqref{cond1} and \eqref{as3})
\begin{equation}\label{mem1}
\begin{split}
&{}^{\widetilde{h}}\mathbf{Ric}_{ab}=\frac{1}{2\widetilde{X}^2}(\nabla_a \widetilde{X}\nabla_b \widetilde{X}+\nabla_a \widetilde{Y}\nabla_b \widetilde{Y}),\\
&{}^{\widetilde{h}}\square (\widetilde{X}+i\widetilde{Y})=\frac{1}{\widetilde{X}}\widetilde{h}^{ab}\partial_a(\widetilde{X}+i\widetilde{Y}) \partial_b(\widetilde{X}+i\widetilde{Y}),\\
&\widetilde{X}^2(\widetilde{\nabla}_a\widetilde{A}_b-\widetilde{\nabla}_b\widetilde{A}_a)=\widetilde{\in}_{abc}\widetilde{\nabla}^c\widetilde{Y},
\end{split}
\end{equation}
in $\mathcal{D}$. Then we construct the space-time metric $\widetilde{\g}$ in $\mathcal{D}\times I$ as in \eqref{spa1}--\eqref{spa3}. In view of Proposition \ref{mainprop2} (a), it remains to show that we can arrange our construction in such a way that the vector-field $Z$ cannot be extended as a Killing vector-field for the modified metric $\widetilde{\g}$. Using Proposition \ref{mainprop2} (b), it suffices to prove that we can arrange the construction in subsection \ref{mem6} such that the vector-field $\partial_3$ cannot be extended to a vector-field $Z'$ in $\mathcal{D}$ such that
\begin{equation}\label{mem2}
\Lie_{Z'}\widetilde{h}=0\,\,\text{ and }\,\,Z'(\widetilde{X})=Z'(\widetilde{Y})=0\,\,\text{ in }\mathcal{D}.
\end{equation}
More precisely, we assume that \eqref{mem2} holds and show that there is a choice of $\widetilde{X}, \widetilde{Y}$ along $\mathcal{N}_1$ such that \eqref{mem1} is violated.

Assuming that \eqref{mem2} holds, we define the geodesic vector-field $\widetilde{L}$ in $\mathcal{D}$ as in subsection \ref{mem6} and notice that
\begin{equation*}
\widetilde{\nabla}_{\widetilde{L}}\widetilde{h}(\widetilde{L},Z')=0.
\end{equation*}
Recall that, see \eqref{L2},
\begin{equation*}
\widetilde{h}(\widetilde{L},\widetilde{L})=0,\qquad [\widetilde{L},Z']=0,\qquad \widetilde{h}(\widetilde{L},Z')=-1,\qquad \text{ along }\mathcal{N}_0.
\end{equation*}
Since $\widetilde{h}(\widetilde{L},Z')=-1$ along $\mathcal{N}_0$, it follows that
\begin{equation*}
\widetilde{h}(\widetilde{L},Z')=-1\qquad\text{ in }\mathcal{D}.
\end{equation*}
We let $e_{(2)}:=\widetilde{L}$, $e_{(3)}:=Z'$, and fix an additional smooth vector-field $e_{(1)}$ in $\mathcal{D}$ such that $\widetilde{h}(e_{(1)},e_{(2)})=\widetilde{h}(e_{(1)},e_{(3)})=\widetilde{h}(e_{(1)},e_{(1)})-1=0$, i.e. 
\beaa
e_{(1)}^a=\widetilde{\in}^{abc} \Lt_b  Z'_c.
\eeaa 
To summarize, assuming \eqref{mem2}, we have constructed a frame $e_{(1)},e_{(2)},e_{(3)}$ in $\mathcal{D}$ such that
\begin{equation}\label{mem4}
\widetilde{h}(e_{(1)},e_{(1)})-1=\widetilde{h}(e_{(1)},e_{(2)})=\widetilde{h}(e_{(1)},e_{(3)})=\widetilde{h}(e_{(2)},e_{(2)})=\widetilde{h}(e_{(2)},e_{(3)})+1=0.
\end{equation}
We define the connection coefficients
\begin{equation*}
\Gamma_{(a)(b)(c)}=\widetilde{h}(e_{(a)},\widetilde{\nabla}_{e_{(c)}}e_{(b)}).
\end{equation*}
Using the identities $\widetilde{\nabla}_{\widetilde{L}}\widetilde{L}=0$ and $\Lie_{Z'}\widetilde{h}=0$, it follows that
\begin{equation*}
\Gamma_{(a)(2)(2)}=0\,\text{ for any }\,a\in\{1,2,3\},\quad\Gamma_{(a)(3)(c)}+\Gamma_{(c)(3)(a)}=0\,\text{ for any }\,(a,c)\in\{1,2,3\}^2.
\end{equation*}
Since $[\widetilde{L},Z']=0$ along $\mathcal{N}_0$ and
\begin{equation*}
0=\Lie_{Z'}(\widetilde{L}^a\widetilde{\nabla}_a\widetilde{L}_b)=\widetilde{L}^a\widetilde{\nabla}_a(\Lie_{Z'}\widetilde{L}_b)+\widetilde{h}^{ac}\widetilde{\nabla}_a\widetilde{L}_b(\Lie_{Z'}\widetilde{L}_c),
\end{equation*}
it follows that $\Lie_{Z'}\widetilde{L}=0$ in $\mathcal{D}$. Then, using the definition of $e_{(1)}$, it follows that $\Lie_{Z'}e_{(1)}=0$ in $\mathcal{D}$, therefore
\begin{equation*}
\Gamma_{(a)(3)(c)}=\Gamma_{(a)(c)(3)}\,\text{ for any }\,(a,c)\in\{1,2,3\}^2.
\end{equation*}
To summarize, letting $F=\widetilde{h}(Z',Z')=\widetilde{h}_{(3)(3)}$, we have
\begin{equation}\label{mem7.5}
\begin{split}
&\widetilde{h}_{(1)(1)}-1=\widetilde{h}_{(1)(2)}=\widetilde{h}_{(1)(3)}=\widetilde{h}_{(2)(2)}=\widetilde{h}_{(2)(3)}+1=0,\quad \widetilde{h}_{(3)(3)}=F,\\
&\widetilde{h}^{(1)(1)}-1=\widetilde{h}^{(1)(2)}=\widetilde{h}^{(1)(3)}=\widetilde{h}^{(3)(3)}=\widetilde{h}^{(2)(3)}+1=0,\quad \widetilde{h}^{(2)(2)}=-F,
\end{split}
\end{equation}
and
\begin{equation}\label{mem8}
\begin{split}
&\Gamma_{(1)(1)(1)}=\Gamma_{(2)(2)(2)}=\Gamma_{(3)(3)(3)}=\Gamma_{(1)(1)(2)}=\Gamma_{(1)(2)(2)}=\Gamma_{(2)(1)(2)}=\Gamma_{(2)(2)(1)}=0,\\
&\Gamma_{(3)(2)(2)}=\Gamma_{(2)(3)(2)}=\Gamma_{(2)(2)(3)}=\Gamma_{(3)(1)(1)}=\Gamma_{(1)(3)(1)}=\Gamma_{(1)(1)(3)}=0,\\
&-\Gamma_{(1)(2)(1)}=\Gamma_{(2)(1)(1)},\quad -\Gamma_{(3)(3)(a)}=-\Gamma_{(3)(a)(3)}=\Gamma_{(a)(3)(3)}=-\frac{1}{2}e_{(a)}(F),\quad a\in\{1,2\},\\
&-\Gamma_{(3)(1)(2)}=\Gamma_{(1)(3)(2)}=-\Gamma_{(2)(3)(1)}=\Gamma_{(3)(2)(1)}=-\Gamma_{(2)(1)(3)}=\Gamma_{(1)(2)(3)},\\
&e_{(3)}(F)=e_{(3)}(\Gamma_{(a)(b)(c)})=0,\qquad [e_{(3)},e_{(a)}]=0,\qquad a,b,c\in\{1,2,3\}.
\end{split}
\end{equation}

We derive now several identities for the connection coefficients $\Gamma$ and the curvature ${}^{\widetilde{h}}\mathbf{R}$. Clearly
\begin{equation*}
\begin{split}
&{}^{\widetilde{h}}\mathbf{R}_{(a)(b)(c)(d)}=\widetilde{h}(e_{(a)},[\widetilde{\nabla}_{e_{(c)}}(\widetilde{\nabla}_{e_{(d)}}e_{(b)})-\widetilde{\nabla}_{e_{(d)}}(\widetilde{\nabla}_{e_{(c)}}e_{(b)})-\widetilde{\nabla}_{[e_{(c)},e_{(d)}]}e_{(b)}])\\
&=\widetilde{h}(e_{(a)},[\widetilde{\nabla}_{e_{(c)}}({\Gamma^{(m)}}_{(b)(d)}e_{(m)})-\widetilde{\nabla}_{e_{(d)}}({\Gamma^{(m)}}_{(b)(c)}e_{(m)})-({\Gamma^{(m)}}_{(d)(c)}-{\Gamma^{(m)}}_{(c)(d)})\widetilde{\nabla}_{e_{(m)}}e_{(b)}])\\
&=e_{(c)}(\Gamma_{(a)(b)(d)})-e_{(d)}(\Gamma_{(a)(b)(c)})\\
&+{\Gamma^{(m)}}_{(b)(d)}\Gamma_{(a)(m)(c)}-{\Gamma^{(m)}}_{(b)(c)}\Gamma_{(a)(m)(d)}+({\Gamma^{(m)}}_{(c)(d)}-{\Gamma^{(m)}}_{(d)(c)})\Gamma_{(a)(b)(m)}
\end{split}
\end{equation*}
for any $a\in\{1,2\}$ and $b,c,d\in\{1,2,3\}$. Using also the identities \eqref{mem7.5} and \eqref{mem8}, it follows that
\begin{equation}\label{mem10}
\begin{split}
&{}^{\widetilde{h}}\mathbf{R}_{(a)(3)(2)(3)}=e_{(2)}(\Gamma_{(a)(3)(3)})-\Gamma_{(2)(3)(3)}\Gamma_{(a)(3)(2)}-\Gamma_{(1)(3)(2)}\Gamma_{(a)(1)(3)}+\Gamma_{(3)(3)(2)}\Gamma_{(a)(2)(3)},\\
&{}^{\widetilde{h}}\mathbf{R}_{(1)(2)(2)(3)}=e_{(2)}(\Gamma_{(1)(2)(3)}),\\
&{}^{\widetilde{h}}\mathbf{R}_{(2)(1)(2)(1)}=e_{(2)}(\Gamma_{(2)(1)(1)})+\Gamma_{(2)(1)(1)}\Gamma_{(1)(2)(1)}.
\end{split}
\end{equation}

We can now obtain our desired  contradiction by   constructing  a pair of smooth functions $\widetilde{X},\widetilde{Y}$ along $\mathcal{N}_1$  such that not all the identities 
 above (starting with  \eqref{mem1}) can  be  simultaneously verified along $\mathcal{N}_1$. For this we fix a smooth system of coordinates $y=(y^1,y^2,y^3)$ in a neighborhood of the point $p$ in $\Pi$ such that
\begin{equation*}
\mathcal{N}_1=\{q:y^3(q)=0\},\qquad \mathcal{N}_0=\{q:y^2(q)=0\},\qquad L=\widetilde{L}=\frac{d}{dy^2}\,\text{ along }\,\mathcal{N}_1.
\end{equation*}
More precisely  we fix the $L$,  as in the unperturbed Kerr, in a neighborhood of $p$ and define 
first $y^2$ such that $y^2$ vanishes on $\NN_0$ and $L(y^2)=1$. Then we   complete the coordinate system  on $\NN_0$ and extend  it by solving $L(y^1)=L(y^3)=0$. 

Assume $\psi:\mathbb{R}^3\to[0,1]$ is a smooth function equal to $1$ in the unit ball and vanishing outside the ball of radius $2$. We are looking for functions $\widetilde{X},\widetilde{Y}$ of the form
\begin{equation}\label{mem12}
\widetilde{X}(q)=X(q),\qquad\widetilde{Y}(q)=Y(q)+\eps\psi((y(q)-y(p'))/\eps),\qquad q\in\mathcal{N}_1,
\end{equation}
where $p'$ is a fixed point in $\mathcal{N}_1\cap D_+$ sufficiently close to $p$, and $(X,Y)$ are as in \eqref{pote}. We show below that such a choice leads to a contradiction, for $\eps$ sufficiently small.

Let
\begin{equation*}
V_1=\frac{d}{dy^1},\qquad V_2=\frac{d}{dy^2},\qquad V_3=\frac{d}{dy^3},\qquad e_{(a)}=K^i_{(a)}V_i.
\end{equation*}
In view of the definitions,
\begin{equation*}
K_{(2)}^1=K_{(2)}^3=K_{(1)}^3=K_{(2)}^2-1=0\qquad\text{ along }\mathcal{N}_1.
\end{equation*}

We use now the last identity  in \eqref{mem10} and the first identity in \eqref{mem1}, along $\mathcal{N}_1$. Since ${}^{\widetilde{h}}\mathbf{R}_{(2)(1)(2)(1)}={}^{\widetilde{h}}\mathbf{Ric}_{(2)(2)}$,
 and recalling \eqref{mem7.5} and \eqref{mem8},  we derive
\begin{equation}\label{mem15}
\begin{split}
&V_2(\Gamma_{(2)(1)(1)})-(\Gamma_{(2)(1)(1)})^2=\frac{1}{2\widetilde{X}^2}[V_2(\widetilde{X})^2+V_2(\widetilde{Y})^2],
\end{split}
\end{equation}
along $\mathcal{N}_1$. In addition, since 
\begin{equation}
[e_{(2)},e_{(1)}]=[V_2,K^1_{(1)}V_1+K^2_{(1)}V_2]=V_2(K^1_{(1)})V_1+V_2(K^2_{(1)})V_2
\label{eq:mem15}
\end{equation}
along $\mathcal{N}_1$, it follows that
\begin{equation}\label{mem16}
\begin{split}
&V_2(K^1_{(1)})=K^1_{(1)}\cdot\widetilde{h}([e_{(2)},e_{(1)}],e_{(1)})=K^1_{(1)}\Gamma_{(2)(1)(1)},
\end{split}
\end{equation}
along $\mathcal{N}_1$.
 Using the ansatz  \eqref{mem12} together with  \eqref{mem15}, and \eqref{mem16}, it follows that
\begin{equation}\label{mem17}
|G|+|V_2(G)|\lesssim 1\text{ for any }G\in \{\Gamma_{(2)(1)(1)},K_{(1)}^1,1/K_{(1)}^1\},
\end{equation}
along $\mathcal{N}_1$, uniformly for all $p'\in\mathcal{N}_1$ sufficiently close to $p$ and $\eps\leq \eps(p')$ sufficiently small.

 Next we use  the identity on  the second line of 
  \eqref{mem10}   and the  Ricci  identity in \eqref{mem1}, along $\mathcal{N}_1$. Since  ${}^{\widetilde{h}}\mathbf{R}_{(1)(2)(2)(3)}=-{}^{\widetilde{h}}\mathbf{Ric}_{(1)(2)}$,  and  recalling \eqref{mem7.5},  \eqref{mem8}  we infer that,
\begin{equation}\label{mem15.1}
V_2(\Gamma_{(1)(2)(3)})=\frac{-1}{2\widetilde{X}^2}[V_2(\widetilde{X})\cdot (K_{(1)}^1V_1+K^2_{(1)}V_2)(\widetilde{X})+V_2(\widetilde{Y})\cdot (K_{(1)}^1V_1+K^2_{(1)}V_2)(\widetilde{Y})],
\end{equation}
along $\mathcal{N}_1$. In addition, using again \eqref{eq:mem15},   it follows that
\begin{equation}\label{mem16.1}
V_2(K^2_{(1)})=-\widetilde{h}([e_{(2)},e_{(1)}],e_{(3)})+V_2(K_{(1)}^1)\widetilde{h}(V_1,e_{(3)})=2\Gamma_{(1)(2)(3)}+K_{(1)}^2V_2(K_{(1)}^1)/K_{(1)}^1
\end{equation}
along $\mathcal{N}_1$. Using once more the ansatz   \eqref{mem12} together with   \eqref{mem15.1}, and \eqref{mem16.1} as well the previously  established bounds \eqref{mem17}, it follows that
\begin{equation}\label{mem17.1}
|G|+|V_2(G)|\lesssim 1\text{ for any }G\in \{\Gamma_{(2)(1)(1)},K_{(1)}^1,1/K_{(1)}^1,\Gamma_{(1)(2)(3)},K_{(1)}^2\},
\end{equation}
along $\mathcal{N}_1$, uniformly for all $p'\in\mathcal{N}_1$ sufficiently close to $p$ and $\eps\leq \eps(p')$ sufficiently small. 

Using the  Ricci identity in \eqref{mem1}, the identities $e_{(3)}(\widetilde{X})=e_{(3)}(\widetilde{Y})=0$, and the bounds \eqref{mem17.1}, it follows that
\begin{equation*}
\sum_{a,b\in\{1,2,3\}}|{}^{\widetilde{h}}\mathbf{Ric}_{(a)(b)}|\lesssim 1\quad\text{ along }\mathcal{N}_1.
\end{equation*}
Using   the first  identity in \eqref{mem10} with $a=2$, the identity 
\begin{equation*}
{}^{\widetilde{h}}\mathbf{R}_{(2)(3)(2)(3)}={}^{\widetilde{h}}\mathbf{Ric}_{(2)(3)}+(1/2)({}^{\widetilde{h}}\mathbf{Ric}_{(1)(1)}+F{}^{\widetilde{h}}\mathbf{Ric}_{(2)(2)}),
\end{equation*}
and \eqref{mem8}, it follows that
\begin{equation*}
\begin{split}
&V_2(F)=-2\Gamma_{(2)(3)(3)},\\
&V_2(\Gamma_{(2)(3)(3)})=-(\Gamma_{(1)(2)(3)})^2+{}^{\widetilde{h}}\mathbf{Ric}_{(2)(3)}+(1/2)({}^{\widetilde{h}}\mathbf{Ric}_{(1)(1)}+F\cdot{}^{\widetilde{h}}\mathbf{Ric}_{(2)(2)}).
\end{split}
\end{equation*}
Using again \eqref{mem17.1}, it follows that
\begin{equation}\label{mem20}
|F|+|V_2(F)|+|V_2(V_2(F))|\lesssim 1\quad\text{ along }\mathcal{N}_1,
\end{equation}
uniformly for all $p'\in\mathcal{N}_1$ sufficiently close to $p$ and $\eps\leq \eps(p')$ sufficiently small. 

We can now derive a contradiction by examining the second equation in \eqref{mem1},
\begin{equation*}
\widetilde{h}^{(a)(b)}\widetilde{\nabla}_{(a)}\widetilde{\nabla}_{(b)}(\widetilde{Y})=2\widetilde{X}^{-1}\widetilde{h}^{(a)(b)}e_{(a)}(\widetilde{X})e_{(b)}(\widetilde{Y}).
\end{equation*}
Using \eqref{mem17} and \eqref{mem20}, it follows that
\begin{equation*}
|e_{(1)}(e_{(1)}(\widetilde{Y}))-Fe_{(2)}(e_{(2)}(\widetilde{Y}))|\lesssim 1\quad\text{ along }\mathcal{N}_1,
\end{equation*}
uniformly for all $p'\in\mathcal{N}_1$ sufficiently close to $p$ and $\eps\leq \eps(p')$ sufficiently small.
This cannot happen,  as can easily be seen by letting first  $\ep\to 0$ and then $p'\to p$, taking into account that $F$ and $K_{(1)} ^2$ vanish along $\NN_0\cap \NN_1$.
\end{proof}

\section{Extension across null hypersurfaces}\label{further}

Assume in this section that $(\M, \mathbf{g})$ is a 4-dimensional Lorentzian manifold satisfying the Einstein-vacuum equations $\mathrm{\bf Ric}(\mathbf{g})=0$, $p\in \M$ is a fixed point along a smooth  null hypersurface
 $\underline{\mathcal{N}}\subseteq \M$ (given by the level hypersurface
 of a smooth function $\ub:\M\to \RRR$)  with fixed null vector-field  $\Lb$ at $p$.
  Assume that $u:\M\to\mathbb{R}$ is a smooth optical  function transversal to $\NNb$,  more precisely,
\begin{equation}\label{assum1}
\D^\al u\D_\al u=0\text{ in }\M,\,\,\,   u(p)=0,\,\,\,   (\D^\al\ub\D_\al u)(p)=-1.
\end{equation}
Let  $\NN$ be   the null hypersurface   passing through
$p$ generated by the zero level set of u, i.e.
$
\NN =\{x\in \M/ u(x)=0\}
$
and
$
L=-\g^{\al\be}\partial_\al u\partial_\be$
its  null geodesic generator. Let
\beaa
O_-:=\{x\in \M:\um(x)<0\}
\eeaa
and  assume that $Z$ is a smooth Killing vector-field in $O_-$.

\subsection{An extendibility criterion} We extend $Z$ to neighborhood of $p$ as in \eqref{pr1.1}, such that
\begin{equation*}
L^\al L^\be(\D_\al\D_\be Z_\mu-Z^\rho\R_{\rho\al\be\mu})=0.
\end{equation*}

\begin{theorem}\label{suff3}
Assume that we have, along the null hypersurface  $\NN$,
\begin{equation}\label{suff4}
(\Lie_Z\R)(L,X,L,Y)=0
\end{equation}
for any vector-fields $X,Y\in T(\M)$ tangent to $\mathcal{N}$. Then there is a  neighborhood $U$ of $p$  such that $\Lie_Z\g=0$ in $U$.
\end{theorem}

\begin{remark}\label{suff3.2}
The sufficient condition \eqref{suff4} may be replaced by a sufficient condition at the level of the deformation tensor $\pi$, namely
\begin{equation}\label{suff4.2}
(\Lie_Z\g)(X,Y)=0\qquad \text{ along }\NN,
\end{equation}
for any vector-fields $X,Y\in T(\M)$ tangent to $\mathcal{N}$. Both \eqref{suff4} and \eqref{suff4.2} lead to the conclusion \eqref{tro2}, using the identities \eqref{tro100}--\eqref{tr104}.
\end{remark}

\begin{proof}[Proof of Theorem \ref{suff3}]
 According to  the results proved in the section  \ref{sec:tens-eq}  we introduce the tensors $W$, $\pi$, $\omega$, $B$, $\dot{B}$ and $P$ as in Definition \ref{pro0}. Recall that, see Lemma \ref{pro3},
\begin{equation}\label{tro1}
 \pi_{\a\mu}L^\mu=0,  \quad   \omega_{\a\mu}L^\mu=0, \quad    P_{\a\b\mu} L^\mu=0.
\end{equation}
Since $B=\frac 1 2 (\pi+\omega)$ we also have $ B_{\a\mu}L^\mu=0$.
We fix a function $y:\mathcal{N}\to\mathbb{R}$ such that $y$ vanishes on $\mathcal{N}\cap\underline{\mathcal{N}}$ and $L(y)=1$ along $\mathcal{N}$. Then we fix a frame $(e_1,e_2,e_3,e_4)$ along $\mathcal{N}$ such that
\begin{equation*}
\begin{split}
&e_1,e_2,e_4\text{ are tangent to }\mathcal{N},\qquad e_4=L,\qquad e_1(y)=e_2(y)=0,\\
&\g(e_1,e_2)=\g(e_a,e_a)-1=\g(e_4,e_3)+1=\g(e_a,e_3)=\g(e_3,e_3)=0,\,\,\, a\in\{1,2\}.
\end{split}
\end{equation*}

Our main goal is to show that the tensors $W,B,\dot{B},P$ vanish along $\mathcal{N}$. For any tensor $M=M_{\al_1.\ldots\al_k}$ and any $s\in\mathbb{Z}$ we define $M^{\geq s}$ any component of the tensor $M$ in the basis $(e_1,e_2,e_3,e_4)$ of signature $\geq s$, where the signature of the component $M_{\al_1\ldots\al_k}$ is equal to the difference between the number of $4$'s and the number of $3$'s in $(\al_1,\ldots\al_k)$. Thus, for example, $\dot{B}^{\geq 0}\in\{\dot{B}_{44},\dot{B}_{4a},\dot{B}_{a4},\dot{B}_{43},\dot{B}_{34},\dot{B}_{ab}:a,b\in\{1,2\}\}$.

Recall our main transport equations, see Lemma \ref{pro5},
\begin{equation}\label{tro100}
\D_LB_{\al\be}=\dot{B}_{\al\be},
\end{equation}
\begin{equation}\label{tro101}
   \D_L\dot{B}_{\a\b}=L^\mu L^\nu   (\Lie_Z \R)_{\mu\al\be\nu}-2\dot{B}_{\nu\be}\D_\al L^\nu-{\pi_\be}^\rho L^\mu L^\nu\R_{\mu\al\rho\nu},
      \end{equation}
and
\begin{equation}\label{tro102}
\D_L P_{\al\be\mu}=2L^\nu W_{\al\be\mu\nu}+2L^\nu {B_{\mu}}^\rho\R_{\al\be\rho\nu}-\D_\mu L^\rho P_{\al\be\rho},
\end{equation}
and our main divergence equation, see Lemma \ref{pro2},
\begin{equation}\label{tr103}
\begin{split}
\D^\a W_{\a\b\ga\de}=\frac{1}{2}\big(&B^{\mu\nu}\D_\nu \R_{\mu\b\ga\de}+\g^{\mu\nu}P_{\mu\rho\nu}\R^\rho\,_{\b\ga\de}\\
  &+P_{\b\nu\mu}\R^{\mu\nu}\,_{\ga\de}+P_{\ga\nu\mu}   \R^{\mu}\,_\b\,^\nu\,_{\de}
  +P_{\de\nu\mu}\R^{\mu}\,_{\b\ga}\,^\nu\big).
\end{split}
\end{equation}
In view of the definitions we also have
\begin{equation}\label{tr104}
\D_\al L_4=0,\qquad \D_4 L_\al=0,\qquad \al\in\{1,2,3,4\}.
\end{equation}

We use equations \ref{tro100}, \eqref{tro101}, and \eqref{tr104}, together with the assumption $\Lie_Z \R_{4ab 4}=0$, $a,b\in\{1,2\}$ to write, schematically,
\begin{equation*}
(\D_L B)^{\geq 0}=\mathcal{M}(\dot{B}^{\geq 0}),\qquad (\D_L \dot{B})^{\geq 0}=\mathcal{M}(\dot{B}^{\geq 0})+\mathcal{M}(B^{\geq 0}).
\end{equation*}
Therefore
\begin{equation*}
B^{\geq 0}=0,\qquad \dot{B}^{\geq 0}=0.
\end{equation*}
According to  Proposition \ref{pro1},
$W=\Lie_Z \R-\frac 1 2  B\odot \R$,
  using again  the assumption \eqref{suff4}, it follows that
\begin{equation*}
W^{\geq 2}=0.
\end{equation*}
Using now \eqref{tro102} and  the identity $P_{\al\be4}=0$, see \eqref{tro1}, it follows that $P^{\geq 1}=0$. Therefore
\begin{equation}\label{tro2}
B^{\geq 0}=0,\qquad\dot{B}^{\geq 0}=0,\qquad P^{\geq 1}=0,\qquad W^{\geq 2}=0\qquad \text{ along }\NN.
\end{equation}

Using \eqref{tro2} and the general symmetries of Weyl fields, equation \eqref{tr103} with $(\be\ga\de)=(4a4)$, $a\in\{1,2\}$, gives, schematically,
\begin{equation*}
(\D_L W)^{\geq 1}=\mathcal{M}(B^{\geq {-1}})+\mathcal{M}(P^{\geq 0})+\mathcal{M}(W^{\geq 1}).
\end{equation*}
Using the transport equations \eqref{tro100}, \eqref{tro101}, and \eqref{tro102}, together with the identities \eqref{tr104} and \eqref{tro2} we derive, schematically,
\begin{equation*}
\begin{split}
&(\D_LB)^{\geq -1}=\mathcal{M}(\dot{B}^{\geq -1}),\\
&(\D_L\dot{B})^{\geq -1}=\mathcal{M}(W^{\geq 1})+\mathcal{M}(\dot{B}^{\geq -1})+\mathcal{M}(B^{\geq -1}),\\
&(\D_LP)^{\geq 0}=\mathcal{M}(W^{\geq 1})+\mathcal{M}(B^{\geq -1})+\mathcal{M}(P^{\geq 0}).
\end{split}
\end{equation*}
Therefore, \eqref{tro2} can be upgraded to
\begin{equation}\label{tro3}
B^{\geq -1}=0,\qquad\dot{B}^{\geq -1}=0,\qquad P^{\geq 0}=0,\qquad W^{\geq 1}=0\qquad \text{ along }\NN.
\end{equation}

We can now continue this procedure. Using \eqref{tro3} and the general symmetries of Weyl fields, equation \eqref{tr103} with $(\be\ga\de)=(434)$ and $(\be\ga\de)=(412)$ gives, schematically,
\begin{equation*}
(\D_L W)^{\geq 0}=\mathcal{M}(B^{\geq {-2}})+\mathcal{M}(P^{\geq -1})+\mathcal{M}(W^{\geq 0}).
\end{equation*}
The transport equations \eqref{tro100}, \eqref{tro101}, and \eqref{tro102}, together with the identities \eqref{tr104} and \eqref{tro3} give, schematically,
\begin{equation*}
\begin{split}
&(\D_LB)^{\geq -2}=\mathcal{M}(\dot{B}^{\geq -2}),\\
&(\D_L\dot{B})^{\geq -2}=\mathcal{M}(W^{\geq 0})+\mathcal{M}(\dot{B}^{\geq -2})+\mathcal{M}(B^{\geq -2}),\\
&(\D_LP)^{\geq -1}=\mathcal{M}(W^{\geq 0})+\mathcal{M}(B^{\geq -2})+\mathcal{M}(P^{\geq -1}).
\end{split}
\end{equation*}
Therefore, \eqref{tro3} can be upgraded to
\begin{equation}\label{tro4}
B=0,\qquad\dot{B}=0,\qquad P^{\geq -1}=0,\qquad W^{\geq 0}=0\qquad \text{ along }\NN.
\end{equation}

Using \eqref{tro4} and the general symmetries of Weyl fields, equation \eqref{tr103} with $(\be\ga\de)=(4a3)$, $a\in\{1,2\}$ gives, schematically,
\begin{equation*}
(\D_L W)^{\geq -1}=\mathcal{M}(P^{\geq -2})+\mathcal{M}(W^{\geq -1}).
\end{equation*}
The transport equation \eqref{tro102}, and the identities \eqref{tr104} and \eqref{tro4} gives, schematically,
\begin{equation*}
(\D_LP)^{\geq -2}=\mathcal{M}(W^{\geq -1})+\mathcal{M}(P^{\geq -2}).
\end{equation*}
Therefore, \eqref{tro4} can be upgraded to
\begin{equation}\label{tro5}
B=0,\qquad\dot{B}=0,\qquad P^{\geq -2}=0,\qquad W^{\geq -1}=0\qquad \text{ along }\NN.
\end{equation}

Using \eqref{tro5}, \eqref{tr103} and the general symmetries of Weyl fields, it $\D_3W_{4343}=0$ and $\D_3W_{4312}=0$. Thus $\D_3W_{4a3b}=0$ along $\NN$, $a,b\in\{1,2\}$. Therefore, the divergence equation \eqref{tr103} with $(\be\ga\de)=(a3b)$, $a,b\in\{1,2\}$, and the transport equation \eqref{tro102} give, schematically,
\begin{equation*}
(\D_LW)^{\geq -2}=\mathcal{M}(P^{\geq -3})+\mathcal{M}(W^{\geq -2}),\qquad (\D_LP)^{\geq -3}=\mathcal{M}(P^{\geq -3}).
\end{equation*}
Therefore we proved that
\begin{equation}\label{tro6}
B=0,\qquad\dot{B}=0,\qquad P=0,\qquad W=0\qquad \text{ along }\NN.
\end{equation}

To prove now that $B,\dot{B},P,W$ vanish in a full neighborhood of the point $p$ we use Proposition \ref{equ1}, Lemma \ref{extendedCarl2} and the observation that, for $\eps_0$ sufficiently small, the functions
\begin{equation*}
f_{\pm}=(\um+\eps_0)(\pm u+\eps_0)
\end{equation*}
are strongly pseudo-convex in a sufficiently small neighborhood of the point $p$. See \cite[Appendix A]{AlIoKl} for more details.
\end{proof}

\subsection{A non-extendible example} In this subsection we provide examples showing that Killing vector-fields do not extend, in general, across null hypersurfaces in space-times satisfying the Einstein-vacuum equations.

\begin{theorem}\label{nj2}
With the notation at the beginning of the section, we  further
assume that $Z(u)=0$ in $O_-$ and that $Z$ does not vanish identically in a neighborhood of $p$ in $O_-$. Then there is a neighborhood $U$ of $p$ diffeomorphic to the open ball $B_1\subseteq\mathbb{R}^4$ and a smooth Lorentz metric $\h$ in $U$ such that $\mathrm{\bf Ric}(\h)=0$ in $U$, $\h=\mathbf{g}$ in $O_-$, but $Z$ does not admit an extension as a smooth Killing vector-field for $\h$ in $U$.
\end{theorem}

In other words, the space-time $(\M,\mathbf{g})$ can be modified in a neighborhood $U$ of $p$, on one side of the null hypersurface $\pr O$, in such a way that the resulting space-time is still smooth and satisfies the Einstein-vacuum equations, but the symmetry $Z$ fails to extend to $U$.

\begin{proof}[Proof of Theorem \ref{nj2}]
 We fix a smooth system of coordinates $\Phi^p:B_1\to B_1(p)$, $\Phi^p(0)=p$, where $B_r=\{x\in\mathbb{R}^4:|x|<r\}$, $r>0$, and $B_1(p)$ is an open neighborhood of $p$ in $O$. Let $\partial_1,\ldots,\partial_4$ denote the induced coordinate vector-fields in $B_1(p)$ and let $B_r(p)=\Phi^p(B_r)$, $r\in(0,1]$. For any smooth function $\phi:B\to\mathbb{C}$, where $B\subseteq B_1(p)$ is an open set, and $j=0,1,\ldots$, we define
\begin{equation*}
|\partial^j\phi(x)|=\sum_{\al_1,\ldots,\al_j=1}^4|\partial_{\al_1}\ldots\partial_{\al_j}\phi(x)|,\qquad x\in B.
\end{equation*}
We assume that
\begin{equation*}
\g_{\al\be}(p)=\mathrm{diag}(-1,1,1,1).
\end{equation*}
and, for some constant $A\geq 1$,
\begin{equation}\label{constA}
\begin{split}
\sup_{x\in B_1(p)}\sum_{j=1}^6\Big[|\partial^ju|+|\partial^j\um|+\sum_{\al,\be=1}^4|\partial^j\g_{\al\be}(x)|\Big]\leq A.
\end{split}
\end{equation}

We will construct the neighborhood
\begin{equation*}
U_p=\{x\in B_{\eps_0}(p):u(x)>-\eps_0^2\}
\end{equation*}
for some constant $\eps_0$ sufficiently small (depending only on the constant $A$ in \eqref{constA}). We define first the hypersurface
\begin{equation*}
\mathcal{N}_{0}=\{x\in B_{\eps_0^{1/2}}(p):u(x)=-\eps_0^2\}.
\end{equation*}
 Recall that $L=-\g^{\al\be}\partial_\al u\partial_\be$ and notice that $L$ is tangent to $\mathcal{N}_{0}$. We introduce smooth coordinates $(y^1,y^2,y^4)$ along the hypersurface $\mathcal{N}_{0}$ in such a way that $y^4=0$ on $\underline{\mathcal{N}}\cap\mathcal{N}_{0}$ and $L=\pr_4$, where $\pr_1, \pr_2, \pr_4$ are the induced    coordinate vector-fields along $\mathcal{N}_{0}$.

We consider smooth symmetric tensors $h$ along $\mathcal{N}_{0}$, such that
it coincides with $\g$  on $\mathcal{N}_{0}\cap O_-$ and,   on both sides of $\NN_0$,
\begin{equation}\label{data}
h(\pr_4,\pr_\a)=0, \qquad\text{ in }\mathcal{N}_{0},\qquad \a\in\{1,2,4\}.
\end{equation}
Thus the only nonvanishing components of $h$ are,
\beaa
 h_{ab}=h(\pr_a, \pr_b),\qquad\text{ in }\mathcal{N}_{0},\qquad a,b\in\{1,2\}.
\eeaa
We would like to apply Rendall's theorem \cite[Theorem 3]{Re} to construct the metric $\h$ in the domain of dependence of $\underline{\mathcal{N}}\cup\mathcal{N}_{0}$, such that $\h=\g$ along $ \,\underline{\mathcal{N}}$ and $\g=h$ along $\mathcal{N}_{0}$. The only restriction is that the symmetric tensor $h$ is arranged such that the resulting metric satisfies the Einstein equation
\begin{equation}\label{restr}
\h^{\al\be}\R(L,\partial_\al,L,\partial_\be)=0\qquad\text{ along }\mathcal{N}_{0},
\end{equation}
with $\R$ the Riemann curvature tensor of $\h$.
Recalling the definition of $\R$  and noting that for a  space-time metric $\h$
which coincides with $h$ on $\NN_0$ we must have $\h^{3a}=\h^{33}=0$ and
$\h^{ab}=h^{ab}$, (i.e. $h^{ac} h_{cb}=\de_{ab}$),
 we deduce,
\bea
 \h^{\a\b }\R(L,\partial_\al,L,\partial_\be)&=& -I+II \label{riemannI+II}
 \eea
 \beaa
I&=& \sum_{a,b\in\{1,2\} }  h^{ab}\h(\D_{\pr_a}  \D_{\pr_4}(\pr_4),\pr_b) \\
II&=&\sum_{a,b\in\{1,2\} } h^{ab}\h(\D_{\pr_4}(\D_{\pr_a} \pr_4),\pr_b)
\eeaa
Thus,
imposing the auxiliary condition\footnote{Writing $\D_{\pr_4} \pr_4=\om \pr_4$
and introducing the  null second fundamental form $\chi_{ab}=\h(\D_{\pr_a}\pr_4, \pr_b)$  of $\NN_0$,   the condition reduces to $\om\c \trch\equiv 0$ along $\NN_0$.},
\begin{equation}\label{covar5}
I\equiv 0\qquad\text{ along } \,   \mathcal{N}_{0},
\end{equation}
 equation \eqref{restr} is equivalent to
\begin{equation}\label{restr2'}
\sum_{a,b\in\{1,2\}}h^{ab}h(\D_{\pr_4}(\D_{\pr_a} \pr_4),\pr_b)=0.
\end{equation}
which can be viewed as a constraint equation for the metric $h$ on $\NN_0$.
Indeed we can  introduce    a  covariant differentiation\footnote{Since the metric $h$ is degenerate on $\mathcal{N}_{0}$, this formula only defines the covariant derivatives $\nabla_{X} Y$ up to a multiple of $L=e_4$.}    along $\NN_0$
compatible with $h$ by the formula,
\bea
\label{covar}
h(\nabla_X Y, Z)&=&\frac 1 2 \big[- Zh(X,Y)+Yh(X, Z)+Xh(Y, Z)\big]
\eea
for $X, Y, Z\in\{\partial_1.\partial_2,\partial_4\}$.
With this definition we  observe that \eqref{restr2'} is equivalent to,
\begin{equation}\label{restr2}
\sum_{a,b\in\{1,2\}}h^{ab}h(\nabla_{\pr_4}  (\nabla_{\pr_a} \pr_4),\pr_b)=0.
\end{equation}
In view of the definition \eqref{covar}, for $a\in\{1,2\}$
\begin{equation*}
\begin{split}
&\nabla_{\pr_a} \pr_4=(1/2)h^{cd}(\pr_4h_{ad}) \pr_c+\mathrm{multiple}(\pr_4),\\
&\nabla_{\pr_4} \pr_a=(1/2)h^{cd}(\pr_4h_{ad}) \pr_c+\mathrm{multiple}(\pr_4).
\end{split}
\end{equation*}
Therefore, the identity \eqref{restr2} is equivalent to
\begin{equation}\label{restr3}
\pr_4(h^{ad}\pr_4h_{ad})+(1/2)h^{ab}h^{cd}\pr_4h_{ad} \pr_4h_{bc}=0.
\end{equation}
Letting
\begin{equation*}
h_{ab}=\phi^2\widehat{h}_{ab},\qquad\det(\widehat{h})=\widehat{h}_{11}\widehat{h}_{22}-\widehat{h}^2_{12}=1,
\end{equation*}
and making the observation $\widehat{h}^{ad}\widetilde{\partial}_4\widehat{h}_{ad}=0$, the identity \eqref{restr3} is equivalent to
\begin{equation}\label{restr4}
\pr_4^2\phi+(1/8)\phi\cdot \widehat{h}^{ab}\widehat{h}^{cd}\pr_4\widehat{h}_{ad} \pr_4\widehat{h}_{bc}=0.
\end{equation}
In other words, we may define $\widehat{h}_{ab}$, $a,b\in\{1,2\}$, as an arbitrary smooth positive definite symmetric tensor along $\mathcal{N}_{0}$, with $\widehat{h}_{11}\widehat{h}_{22}-\widehat{h}^2_{12}=1$ and $\widehat{h}_{ab}=(\g_{11}\g_{22}-\g_{12}^2)^{-1/2}\g(\pr_a,\pr_b)$ in $\mathcal{N}_{0}\cap O_-$. We then define $\phi$ according to the equation \eqref{restr4}, and the full tensor $h=\phi^2\widehat{h}$ along $\mathcal{N}_{0}$. Finally we apply Rendall's theorem \cite[Theorem 3]{Re} to construct a smooth space-time metric $\h $ in $\widetilde{U}_p=\{x\in B_{\eps_0}(p):u(x)\geq -\eps_0^2\}$ satisfying the Einstein-vacuum equations and agreeing with $\g$ in $U_p\cap O_-$ and with $h$ along $\mathcal{N}_{0}\cap B_{\eps_0}$.  Since, by construction, the term $II$  vanishes identically on $\NN_0$ it also follows that the metric $\h$ verifies the auxiliary assumption
\eqref{covar5}.  We now interpret condition \eqref{covar5} using the null second fundamental form of  $\NN_0$ with respect to the $\h$ metric,
\bea
\label{nullsec-f}
\chi(X, Y):=\h(\D_X L, Y), \qquad \forall X, Y \quad \text{tangent to }\, \NN_0
\eea
Clearly $\D_L L= \om L$ along $\NN_0$ for some smooth function $\om$. Thus,
\beaa
I&=& \sum_{a,b\in\{1,2\} }  h^{ab}\h(\D_{\pr_a}  \D_{\pr_4}(\pr_4),\pr_b)=
 \sum_{a,b\in\{1,2\} }  \om  h^{ab}\h(\D_{\pr_a} \pr_4,\pr_b)=\om  h^{ab}\chi_{ab}=\om \trch.
\eeaa
Thus \eqref{covar5} takes the form,
\bea
\label{covar55}
\om\c  \trch=0.
\eea
from which we infer that $\om$ must  vanish in $U_p\cap \NN_0$  (i.e. $\D_L L=0$)
if $\trch$ vanishes at most on a set with empty interior in  $(U_p\cap \NN_0)\setminus O_-$.

On the other hand,
\beaa
\chi_{ab}&=&\frac 1 2 \pr_4 h_{ab}=\frac 1 2 \pr_4(\phi^2 \widehat{h}_{ab})=\phi\pr_4\phi \, \widehat{h}_{ab}+\frac 1 2 \phi^2\pr_4 \widehat{h}_{ab}
\eeaa
from which,
\beaa
\trch&=&\phi^{-2}\widehat{h}^{ab}\big(\phi\pr_4\phi \, \widehat{h}_{ab}+\frac 1 2 \phi^2\pr_4 \widehat{h}_{ab}\big)=2 \phi^{-1}\pr_4 \phi.
\eeaa
Also the traceless part  of $\chi$,
\beaa
\chih_{ab}&=&\chi_{ab}-\frac 1 2 \trch h_{ab}=\frac 1 2\phi^2\pr_4 \widehat{h}_{ab}.
\eeaa
Thus equation \eqref{restr4} takes the well known form
\bea
\pr_4\trch+\frac 1 2 (\trch)^2=-|\chih|_h^2
\eea
from which we infer that $\trch$ can only    vanish  in a set with empty interior in    $(U_p\cap\NN_0)\setminus O_-$  if the same holds true for $\chih$.  Thus we can easily
 choose non-trivial  data on $\NN_0$ such that, for our original choice of $L=\pr_4$,
 we have,
 \bea
 \label{covar555}
 \D_L L=0 \qquad \text{in }\, \, \NN_0\cap U_p.
 \eea

It remains to prove that we can arrange $\widehat{h}_{ab}$ on $\mathcal{N}_{0}$ such that $Z$ does not admit an extension to $U_p$ as a Killing vector-field for $\h$. We extend  first the smooth vector-field $L$  from $U_p\cap O_-$ to    all of  $U_p$ such that, consistent with \eqref{covar555},
\begin{equation*}
\D_L L=0\text{ in }   U_p
\end{equation*}
with $\D$ the covariant differentiation associated to the metric $\h$.

Since $Z(u)=0$ in $O_-$ it follows that $[Z,L]=0$ in $O_-$. Assume, for contradiction, that $Z$ admits an extension  to $U_p$ as a Killing vector-field for $\h$. Then, letting $V=\Lie_Z  L$ we compute in $U_p$
\begin{equation*}
L^\rho \D_\rho V_\al= L^\rho\Lie_Z\D_\rho L_\alpha= -V^\rho D_\rho L_\al.
\end{equation*}
Since $V$ vanishes in $U_p\cap O_-$, it  must vanish in all of $U_p$, i.e.
\begin{equation*}
[L, Z]=0\qquad\text{ in }U_p.
\end{equation*}
  In addition, since
\begin{equation*}
L\h(Z, L)=0,
\end{equation*}
we infer that  $Z$ must remain  tangent to the hypersurface $\mathcal{N}_{-\eps_0^2}$. To summarize,  by contradiction,   we have constructed a vector-field $Z$ in $U_p$  tangent  the hypersurface $\mathcal{N}_{0}$ such that, on $\NN_0\cap U_p$,
\begin{equation}\label{rovar1}
\Lie_{Z}h=0,\qquad [L, Z]=0.
\end{equation}

On the other hand, writing $ Z=Z^1\pr_1+Z^2\pr_2+Z^4\pr_4$ in the system of coordinates along $\mathcal{N}_{0}$ introduced before, the
 identity $\Lie_ Zh=0$ in \eqref{rovar1} gives
\begin{equation*}
0= Z(h_{ab})+\pr_aZ^\rho h_{\rho b}+\pr_bZ^\rho h_{a\rho},\qquad a,b\in\{1,2\}.
\end{equation*}
Therefore
\begin{equation*}
 Z(\mathrm{det} (h))=Z(h_{11}h_{22}-h_{12}^2)= -2(\pr_1Z^1+\pr_2Z^2)\mathrm{det}(h).
\end{equation*}
Since $h=(\mathrm{det}h)^{1/2}\widehat{h}$, the identity $\Lie_{Z}h=0$ shows that
\begin{equation}\label{rovar2}
\Lie_{Z}\widehat{h}=(\pr_1Z^1+\pr_2Z^2)\widehat{h}.
\end{equation}
Notice also that $Z$ does not depend on the choice of the tensor $\widehat{h}$, indeed $Z$ is  defined  simply  by  the relation $[L,Z]=0$ in \eqref{rovar1}). Therefore we obtain a contradiction, by choosing $\widehat{h}$ such that \eqref{rovar2} fails at some point in $\mathcal{N}_{0}\setminus O_-$. This completes the proof.
\end{proof}

\end{document}